\documentclass[11pt]{article}
  \def\ZR{{\mathbb R}}

  \def\ZN{{\mathbb N}}

\def\beq{\begin{equation}}
\def\eeq{\end{equation}}
\def\be{\begin{equation}}
\def\ee{\end{equation}}
\def\beqar{\begin{eqnarray}}
\def\eeqar{\end{eqnarray}}
\def\ber{\begin{eqnarray}}
\def\eer{\end{eqnarray}}
\def\berb{\begin{eqnarray*}}
\def\eerb{\end{eqnarray*}}

\def\SO{\mathop{\rm SO}\nolimits}

\def\norm#1.#2.{\|#1\|_{#2}}
\def\Norm#1.#2.{\big\|#1\big\|_{#2}}
\def\NOrm#1.#2.{\bigg\|#1\bigg\|_{#2}}
\def\NORm#1.#2.{\Big\|#1\Big\|_{#2}}
\def\NORM#1.#2.{\Bigg\|#1\Bigg\|_{#2}}

\newcommand{\eps}{\varepsilon}

\def\vec#1{{\mathchoice{\mbox{\boldmath$\displaystyle#1$}}
{\mbox{\boldmath$\textstyle#1$}}
{\mbox{\boldmath$\scriptstyle#1$}}
{\mbox{\boldmath$\scriptscriptstyle#1$}}}}

\def \0b{{\hbox{\boldmath $0$}}}

\newcommand{\ab}{\vec{a}} 

\newcommand{\ccb}{\vec{c}}
\newcommand{\db}{\vec{d}}
\newcommand{\eb}{\vec{e}} \newcommand{\fb}{\vec{f}}
\newcommand{\gb}{\vec{g}} 
\newcommand{\ib}{\vec{i}} 
 \newcommand{\lb}{\vec{l}}

\newcommand{\ssb}{\vec{s}} \newcommand{\tb}{\vec{t}}
 \newcommand{\vb}{\vec{v}}
 
\newcommand{\yb}{\vec{y}} \newcommand{\zb}{\vec{z}}

\newcommand{\Abb}{{\bf A}} \newcommand{\Bbbb}{{\bf B}}
 
 \newcommand{\Fbb}{{\bf F}}
\newcommand{\Gbb}{{\bf G}} \newcommand{\Hbb}{{\bf H}}
\newcommand{\Ibb}{{\bf I}}

\newcommand{\Qbb}{{\bf Q}} \newcommand{\Rbb}{{\bf R}}

 \newcommand{\Pb}{\vec{P}}
\newcommand{\Qb}{\vec{Q}}

\def \alb{\vec{\alpha}} \def \betab{\vec{\beta}}

\newcommand{\cA}{{\cal A}}

 \newcommand{\cR}{{\cal R}}

\def \ftb{\vec{\tilde{f}}}

\def \ptb{\vec{\tilde{p}}}
\def \qtb{\vec{\tilde{q}}}

\def \stb{\vec{\tilde{s}}}

\def \Ftb{\vec{\tilde{F}}}

\newcounter{primjer}[section]
\newcounter{tvrdnja}[section]
\newcounter{zadatak}[section]
\newcommand{\oRbb}{{\bf \overline{R}}}

\usepackage{esint}
\usepackage{amssymb}
\usepackage{graphicx}
\newcommand{\dist}{\mathop {\mbox{\rm dist}}\nolimits}

\newcommand{\tJ}{\tilde{J}}
\newtheorem{theorem}{Theorem}[section]
\newtheorem{corollary}[theorem]{Corollary}

\newtheorem{proposition}[theorem]{Proposition}
\newtheorem{lemma}[theorem]{Lemma}
\newtheorem{remark}[theorem]{Remark}

\let\EEij=\relax 
\let\EEjj=\relax 
\newcommand{\Reshji}{\relax} 
\newcommand{\Reshjj}{\relax}
\newcommand{\Reshjn}{\relax}
\newcommand{\Qu}{q}

\newcommand{\wtRbb}{{\bf\widetilde{R}}}
\newcommand{\whRbb}{{\bf\widehat{R}}}
\newcommand{\wtyb}{{\vec{\widetilde{y}}}}
\newcommand{\hRbb}{{\bf{R}}^\eps}
\newcommand{\hsb}{{\vec{s}}^\eps}

\newcommand{\Pbh}{{\Pb^{(h)}}}
\newcommand{\Pbhj}{{\Pb^{(h_j)}}}
\newcommand{\Id}{\ib}

\begin{document}

\title{Derivation of the nonlinear bending--torsion model for a junction of elastic rods
}




\author{Josip Tamba\v{c}a         \and
        Igor Vel\v{c}i\'{c} 
}



\maketitle

\begin{abstract}
In this paper we derive the one--dimensional bending--torsion equilibrium model modeling the junction of straight rods. The starting point is a three-dimensional nonlinear elasticity equilibrium problem written as a minimization problem for a union of thin rod--like bodies.
By taking the limit as the thickness of the 3D rods goes to zero, and by using ideas from the theory of $\Gamma$--convergence, we obtain that the resulting model consists of the union of the usual one--dimensional nonlinear bending-torsion rod models which satisfy the following transmission  conditions at the junction point: continuity of displacement and rotation of the cross--sections and balance of contact forces and contact couples.
\end{abstract}

\section{Introduction}

In many real-life situations, such as, for example, in certain types of bridges or building structures, two (or several) elastic rods are connected at one point. Such points where several rods meet are called junctions. Such multiple rods systems may be as small as two rods joining in a non-smooth way, or as complex as several hundreds of interconnected rods forming a massive network. In either case, the basic principles of analysis are the same (although the complexity of the computation depends on the complexity of the system). Therefore, in the present paper, we limit our study to the case of one junction point.

In this paper we consider the equilibrium problem of a three-dimensional elastic body which consists of $n$ straight thin rod-like bodies connected in a single point. Since the rods are thin, the behavior of each rod should be well approximated by the one--dimensional rod model. In order to close the model one needs conditions at the junction point. These conditions can be seen as transmission conditions as well. Since we are interested in the bending-torsion behavior of rods, such a rod is expected to be governed by the fourth order equation, see \cite{Antman}. Since this equation can be written as a first order system in terms of the contact force, the contact couple, the rotation of the cross-section and the deformation (displacement), we expect the following four junction conditions (based on the continuity of the deformation and equilibrium laws):
\begin{itemize}
\item[1)] the sum of all contact forces at the junction is zero,
\item[2)] the sum of all contact couples at the junction is zero,
\item[3)] continuity of the rotation of the cross-section (the angles at the junction point are preserved),
\item[4)] continuity of the displacement (deformation/position) at the junction point.
\end{itemize}
These conditions follow physical intuition and are already used in modelling networks of elastic rods (see e.g. \cite{Tucker1} and \cite{TKCP}; in the case of strings see \cite{zuazua}).

In \cite{Griso1} these junction conditions have been mathematically justified for the case  when the starting configuration is that of three--dimensional {\sl linearized} elasticity.
We justify these junction conditions starting from a three--dimensional nonlinearly hyperelastic material in the formulation as the energy minimization problem.
Since the rods are thin, we recognize the small parameter $h$ describing their thickness.
The mechanical response of rods strongly depends on the relative
magnitude of the applied load with respect to the rod thickness $h$. In
\cite{Muller1} a bending-torsion model of a single nonlinearly elastic
inextensible rod was derived by the theory of $\Gamma$-convergence and
the geometric rigidity theorem from \cite{Muller0}. In order to obtain the bending--torsion model the main assumptions is that the energy of the rod is of the order $h^4$. For other models see \cite{talijani, Muller4}.
In the present paper we would like to obtain junction conditions at the
junction of rods for the case when the total energy functional is of order $h^4$.
However, unlike in the case of a single rod studied in
\cite{Muller0}, in the case when 2 or more rods meet at a junction, we cannot
rescale our problem in such a way that the entire problem is defined
on a canonical domain independent of $h$, at least in a simple way.  To deal with the
complications related to the geometry at the junction, we assume that
the junction region of the rods forms a domain which scales with $h$
(say a sphere, at which all the rods are connected). Then, as $h\to
0$, the junction region converges to a point. This leads to a problem
with no obvious simple canonical domain, and so the results from \cite{Muller1}
cannot be applied directly to this problem. To get around this
difficulty we adapt the ideas from \cite{Muller1} to this new scenario
and express the asymptotic behavior of minimizers in norms depending on the thickness $h$.

Following \cite{Muller1}, we first prove a compactness result (Theorem Theorem~\ref{glavni}) for the sequence of energy minimizers $\yb^{(h)}$
deriving the asymptotic behavior of $\nabla \yb^{(h)}$. Moreover we prove that the
rotations of the cross-sections need to be continuous at the junction
point in the limit as $h\to 0$. Since we are considering a pure
traction problem for rods joining at a point, we still need to control
the displacement of the entire structure. Under the assumption that
the translation of the whole structure is controlled at the end of one
rod, in Lemma~\ref{lemasranje} and Corollary~\ref{ct3.1l4.1} we derive the asymptotic behavior of the minimizing sequence $\yb^{(h)}$ and we obtain that in the
limit, the displacement (deformation) of the rods at the junction
point is continuous. Finally, in Theorem~\ref{zadnji} we derive the model for
the junction of rods.

Junction of elastic rods has been studied by several authors. However,
most results are restricted to linearized elasticity. The first study
of the junction of two rods is given by Le Dret in \cite{LD0}, see also \cite{LD} and \cite{Tambaca}. For systems of rods see  also \cite{Nazarov} and \cite{Panasenko} and references therein. The junction of two plates is studied in \cite{Griso2} and \cite{LD1, LD2}, while \cite{S-P} deals with the junction of beams and
plates. The case of the junction of a three dimensional domain and a two dimensional one is explored in \cite{Ciarlet2}, see also \cite{Ciarlet0} and references therein.  For the asymptotic analysis of the junction between three-dimensional structures and one-dimensional one see \cite{Kozlov} and \cite{Argatov}.
See also \cite{fluidi1, fluidi2} for the asymptotic analysis of the  problem of junctions of thin pipes filled with a fluid using asymptotic expansion method.

Two efforts in the study of junction problems within nonlinear elasticity are made in \cite{non3P} and \cite{nonPR} using asymptotic expansion method. In \cite{non3P} the model of plate inserted in a three-dimensional elastic body is derived, while in \cite{nonPR} a model of junction of rod and plate is derived.




\section{Setting up the problem}
\setcounter{equation}{0}

The domain of the junction of rods we define as a union of cylinders and the "junction" part. Let $n\in \ZN$ denote the number of rods meeting in junction and let $h>0$. Let each rod be of length $L_i$ with the cross-section $h S_i$, where $S_i \subset \ZR^2$ (open, bounded, connected). Let the junction part is of the form $T^h = h T$, for $T \subset \ZR^3$ open, bounded, connected set. Let $\Qbb_i \in \SO(3)$, $i=1,\ldots, n$. The vector $\tb_i= \Qbb_i \eb_1$ denotes the tangential direction of the $i$-th
rod. Then the domain of the junction of rods is given as
$$
\Omega^h = T^h \cup \bigcup_{i=1}^{n} C^h_i, \qquad C^h_i = \Qbb_i ((h, L_i) \times h S_i).
$$
We assume that the domain $\Omega^h$ is open, bounded, connected and with the Lipschitz boundary. We also assume, as in \cite{Muller1}, for each $i$ that
$$
\int_{S_i} x_2 x_3 dx_2 dx_3=\int_{S_i} x_2 dx_2 dx_3=\int_{S_i} x_3 dx_2 dx_3=0.
$$


Every function $\yb \in W^{1,2} ((a,b);\ZR^3)$ we naturally interpret as an element of $W^{1,2} ((a,b)\times \ZR^2;\ZR^3)$. We also define the mapping $\Pbh : (a,b)\times S_i \to (a,b)\times h S_i$ by $\Pbh(x_1,x_2,x_3)=(x_1,h x_2,h x_3)$ and use it to change between thin and thick domain.


The starting point of our analysis is the equilibrium problem of the junction of rods, i.e. elastic body $\overline{\Omega^h}$. The internal energy of the junction of rods  is given by
$$
E^{(h)}(\yb)= \int_{\Omega^h} W(\nabla \yb(x))dx,
$$
for a deformation $\yb \in W^{1,2} (\Omega^h, \ZR^3)$,
where $W: \mathbb{M}^{3 \times 3} \to[0,+\infty]$ is an internal energy density function. For $W$, as in \cite{Muller1}, is supposed to satisfy
\begin{itemize}
\item $W \in C^0(\mathbb{M}^{3 \times 3})$, $W$ is of class $C^2$ in a neighborhood of $\SO(3)$;
\item $W$ is frame-indifferent, i.e., $W(\Fbb)=W(\Rbb\Fbb)$ for every $\Fbb \in \mathbb{M}^{3 \times 3}$ and $\Rbb \in \SO(3)$;
\item $W(\Fbb) \geq C_W \dist^2(\Fbb,\SO(3)), \ W(\Fbb)=0$ if $\Fbb \in \SO(3)$.
\end{itemize}
We are looking for the one-dimensional bending--torsion model of junction of rods. Thus, motivated by \cite{Muller1}, we will assume that the energy $E^{(h)}$ behaves as $h^4$. Then we analyze the behavior of $E^{(h)}(\yb)/h^4$ and derive the one-dimensional model. This is in \cite{Muller1} obtained by $\Gamma$--convergence, but in the junction problem  there is no obvious and simple canonical domain, a domain that is independent of the thickness $h$. Still, using the ideas and techniques of $\Gamma$--convergence we are able to give the asymptotics (in the form (\ref{kvg})) of the infimizing sequence of the total energy functional and the total energy functional itself.

We shall need the following theorem which can be found in \cite{Muller0}.
\begin{theorem}[on geometric rigidity]\label{tgr}
Let $U \subset \ZR^m$ be a bounded Lipschitz domain, $m \geq 2$. Then there exists a constant $C(U)$ with the following property: for every $\vb\in W^{1,2}(U;\ZR^m)$ there is associated rotation $\Rbb\in \SO(m)$ such that:
\begin{equation}\label{gr}
\|\nabla \vb-\Rbb\|_{L^2(U)} \leq C(U)\|\dist(\nabla \vb,\SO(m)\|_{L^2(U)}.
\end{equation}
\end{theorem}

We will apply this theorem in the next section on subdomains of $\Omega^h$ which are of size $h$ in each direction. This is possible since the constant $C(U)$ in the estimate is independent on the translation and dilatation of $U$.
Let us consider the domain $hU$, for $h>0$. Take $\vb \in W^{1,2}(hU;\ZR^m)$. Then the function $\vb^{(h)}(x) = \frac{1}{h} \vb(h x)$ belongs to $W^{1,2}(U;\ZR^m)$ and satisfies the estimate
$$
\|\nabla \vb^{(h)}-\Rbb\|_{L^2(U)} \leq C(U)\|\dist(\nabla \vb^{(h)},\SO(m)\|_{L^2(U)}.
$$
Since $\nabla \vb^{(h)}=\nabla \vb (hx)$ after the change of variables
in the norms we obtain that the estimate (\ref{gr}) holds for
$\vb$ with the same constant $C(U)$. See also \cite{Muller0}.

Throughout the paper we use the following function space
$$
W^{1,p}(\Omega;\SO(3))=\{ \Rbb \in W^{1,p}(\Omega;\ZR^{3\times 3} | \Rbb(x) \in \SO(3) \ \textrm{for a.e.} \ x \in \Omega \}.
$$
Moreover, by $\| \cdot \|$ (without subscript) we denote the Frobenius matrix norm.

\section{Compactness}
\setcounter{equation}{0}

In this section, following \cite{Muller1} we prove the compactness result (Theorem~\ref{glavni}). Namely,  for $\yb^{(h)}$ that satisfy (\ref{ogranicenost}) (this will be shown for infimizers $\yb^{(h)}$ of the energy of order $h^4$) we obtain asimptotics of $\nabla \yb^{(h)}$. Moreover it turns out that rotations of the cross-sections in the limit, when $h$ tends to 0, need to be continuous in the junction point.
\begin{theorem} \label{glavni}
Let $(\yb^{(h)}) \subset W^{1,2}(\Omega^h;\ZR^3)$ be such that
\begin{equation}\label{ogranicenost}
\limsup_{h \to 0} \frac{1}{h^4} \int_{\Omega^h} \dist^2 (\nabla \yb^{(h)} ,\SO(3))dx < +\infty.
\end{equation}
Then there exists a subsequence (not relabeled) and
$\oRbb_i \in W^{1,2}((0,L_i), \SO(3))$, $i=1,\dots, n$ such that $\oRbb_1(0)=\oRbb_2(0)=\cdots=\oRbb_n (0)$ in the sense of traces and
\begin{equation}\label{3.*}
\lim_{h \to 0} \frac{1}{h^2} \sum_{i=1}^n\int_{C_i^h} \| \nabla \yb^{(h)}(x) - \oRbb_i(x \cdot \tb_i) \|^2 dx = 0.
\end{equation}
\end{theorem}

\begin{prooof}
We follow proof of Theorem 2.2. in \cite{Muller1}.

Now we cover $\Omega^h$ with subdomains of size $h$ in each direction and apply Theorem~\ref{tgr} on each of them.  For every $h>0$ and $i=1,\dots,n$ let $k^h_i \in \mathbb{N}$ be such that $ h \leq L_i/k^h_i <2h $ and let
$$
I^{i}_{a,k^h_i}:=\left( a, a+\frac{L_i}{k^h_i} \right),  \quad a \in [0,L_i) \cap \frac{L_i}{k^h_i}\mathbb{N}.
$$
We apply  Theorem~\ref{tgr} to domains $\Qbb_i((a,a+2h)\times hS_i)$ (when $a=L_i-\frac{L_i}{k^h_i}$ we take $(L_i-2h,L_i)$) and $ T^h \cup \cup_{i=1}^n \Qbb_i((h,2h)\times hS_i)$. Note that $\Qbb_i(I^i_{a,k^k_i}\times hS_i) \subset \Qbb_i((a,a+2h)\times hS_i)$. Then there exist a constant $C$ (independent of $i$ (as there is finite number of domains) and $h$ (by a note after Theorem~\ref{tgr})) and a piecewise constant map $\Rbb^{(h)}:\cup_{i=1}^n \Qbb_i([0,L_i]\times\{0\}\times\{0\})\to \SO(3)$, constant on each $[a, a+\frac{L_i}{k^h_i})$ for $a \in [0,L_i) \cap \frac{L_i}{k^h_i}\mathbb{N}$ and on $T^h \cup \cup_{i=1}^n \Qbb_i([h,\frac{L_i}{k^h_i})\times hS_i)$, such that
for every $i \in \{1, \ldots, n\}$ we have: for  every $a \in [0,L_i) \cap \frac{L_i}{k^h_i}\mathbb{N}$
$$
\int_{\Qbb_i(I^i_{a,k^h_i}\times hS_i)} \|\nabla \yb^{(h)} -\Rbb^{(h)}\|^2dx \leq C \int_{\Qbb_i((a,a+2h)\times hS_i)} \dist^2 (\nabla \yb^{(h)},\SO(3))dx
$$
and
$$
\int_{ T^h \cup \cup_{i=1}^n \Qbb_i((h,\frac{L_i}{k^h_i})\times hS_i)} \hspace{-5ex} \|\nabla \yb^{(h)} -\Rbb^{(h)}\|^2dx \leq C \int_{ T^h \cup \cup_{i=1}^n \Qbb_i((h,2h)\times hS_i)}\hspace{-5ex} \dist^2 (\nabla \yb^{(h)},\SO(3)) dx.
$$
By summing all these estimates, since only neighboring subdomains overlap,  we obtain the inequality
\begin{equation}\label{jjs}
\frac{1}{h^2} \int_{\Omega^h} \|\nabla \yb^{(h)} -\Rbb^{(h)}\|^2dx \leq \frac{2C}{h^2} \int_{\Omega^h} \dist^2 (\nabla \yb^{(h)} ,\SO(3))dx \leq C_1 h^2,
\end{equation}
where the last inequality holds for $h$ small enough by (\ref{ogranicenost}).

In the sequel we show that on a subsequence $\Rbb^{(h)}$ converges to a $W^{1,2}$ function. In order to do that we first estimate the difference of $\Rbb^{(h)}$ on neighboring subdomains.

Let now $a_i \in (0,L_i-4h] \cap \frac{L_i}{k^h_i} \mathbb{N}$,  $b_i=a_i+\frac{L_i}{k^h_i}$. Now we apply Theorem\ref{tgr} on the set $\Qbb_i ((a_i,a_i+4h) \times hS_i)$ we obtain that there exist $\oRbb \in \SO(3)$ such that
$$
\int_{\Qbb_i ((a_i,a_i+4h) \times hS_i)}\|\nabla \yb^{(h)}-\oRbb\|^2 dx \leq C_2 \int_{\Qbb_i ((a_i,a_i+4h) \times hS_i)} \dist^2(\nabla \yb^{(h)},\SO(3))dx.
$$
Then using the facts that $I^i_{a_i,k^h_1},I^i_{b_i,k^h_1}$ are contained in $(a_i,a_i+4h) \times hS_i$ we have for every $i$:
\begin{eqnarray*}
&&\frac{L_i}{k^h_i} \|\Rbb^{(h)}(a_i \tb_i)-\Rbb^{(h)}(b_i \tb_i)\|^2\\
&&\leq 2\frac{L_i}{k^h_i} \left(\|\Rbb^{(h)}(a_i \tb_i)-\Rbb\|^2 +\|\Rbb -\Rbb^{(h)}(b_i \tb_i)\|^2\right)\\
&&\leq 2 \int_{\Qbb_i(I^i_{a_i, h^h_i} \times hS_i)} \|\Rbb^{(h)}(a_i \tb_i)-\oRbb\|^2 + 2\int_{\Qbb_i(I^i_{b_i, h^h_i} \times hS_i)} \|\oRbb -\Rbb^{(h)}(b_i \tb_i)\|^2\\
&&\leq 4 \int_{\Qbb_i(I^i_{a_i, h^h_i} \times hS_i)} \|\Rbb^{(h)}(a_i \tb_i)-\nabla \yb^{(h)}\|^2 + \|\nabla \yb^{(h)}-\oRbb\|^2\\
&& + 4\int_{\Qbb_i(I^i_{b_i, h^h_i} \times hS_i)} \|\oRbb -\nabla\yb^{(h)}\|^2 + \|\nabla\yb^{(h)} -\Rbb^{(h)}(b_i \tb_i)\|^2\\
&&\leq 4 \int_{\Qbb_i(I^i_{a_i, h^h_i} \times hS_i)} \|\Rbb^{(h)}(a_i \tb_i)-\nabla\yb^{(h)}\|^2 + 4 \int_{\Qbb_i ((a_i,a_i+4h) \times hS_i)}\|\nabla \yb^{(h)}-\oRbb\|^2\\
&& + 4\int_{\Qbb_i(I^i_{b_i, h^h_i} \times hS_i)} \|\nabla \yb^{(h)} -\Rbb^{(h)}(b_i \tb_i)\|^2.
\end{eqnarray*}
All terms on the right hand side of the estimate can be estimated by Theorem~\ref{tgr}, so we obtain
\begin{equation} \label{josmalo}
\frac{L_i}{k^h_i} \|\Rbb^{(h)}(a_i \tb_i)-\Rbb^{(h)}(b_i \tb_i)\|^2 \leq \frac{C_3}{h^2} \int_{\Qbb_i((a_i,a_i+4h) \times hS_i)} \dist^2(\nabla \yb^{(h)} ,\SO(3)) dx,
\end{equation}
and similarly, as $I^i_{\frac{L_i}{k^h_i},k^h_i}, T_h$  are contained in $T^h \cup \Qbb_i((h,4h)\times hS_i) $, we obtain
\begin{equation} \label{dosadno5}
\frac{L_i}{k^h_i} \|\Rbb^{(h)}(0)-\Rbb^{(h)}(\frac{L_i}{k^h_i}\tb_i)\|^2 \leq \frac{C_3}{h^2} \int_{T^h \cup \Qbb_i((h,4h) \times hS_i)} \dist^2(\nabla \yb^{(h)} ,\SO(3)) dx.
\end{equation}
Thus we have (since $\Rbb^{(h)}$ is piecewise constant) for every $0 \leq \xi \leq \frac{L_i}{k^h_i}$ and every $i$ and for every $a \in(0,L_i) \cap \frac{L_i}{k^h_i} \mathbb{N}$ s.t $(a,a+4h)\subset (0,L_i)$:
\begin{equation} \label{dosadno1}
\int_{I^i_{a,k^h_1}}\|\Rbb^{(h)}((x_1 +\xi)\tb_i)-\Rbb^{(h)}(x_1 \tb_i)\|^2 dx_1 \leq \frac{C_3}{h^2} \int_{\Qbb_i((a,a+4h)\times hS_i)} \dist^2(\nabla \yb^{(h)},\SO(3))dx,
\end{equation}
since $x_1+\xi$ and $x_1$ belong to the same or neighboring subdomains and we can apply estimate (\ref{josmalo}).
In the same way we can show that for every $i$ and $a$ s.t. $(a-2h,a+2h) \subset (0,L_i)$ and every $-L_i/k^h_i\leq \xi \leq 0$,
\begin{equation}\label{dosadno2}
\int_{I^i_{a,k^h_i}} \|\Rbb^{(h)}((x_1+\xi)\tb_i)-\Rbb^{(h)}(x_1\tb_i)\|^2 dx_1 \leq \frac{C_3}{h^2} \int_{\Qbb_i((a-2h,a+2h)\times hS_i)} \dist^2 (\nabla \yb^{(h)},\SO(3)) dx.
\end{equation}
Let us now look at cylinders $C_1^h$ and $C_2^h$.
By summing estimates (\ref{josmalo}),(\ref{dosadno5}),(\ref{dosadno1}),(\ref{dosadno2}) 
we have
that for every open interval $I'$ compactly contained in $(-L_1,L_2)$ and $\xi \in \ZR$  which satisfies for all $i$, $|\xi| \leq \dist(I',\{-L_1,L_2 \}), |\xi|\leq \frac{L_i}{k^h_i}$
\begin{equation}\label{uzasno1}
\int_{I'} \|\Rbb^{(h)}_m(x_1+\xi)-\Rbb^{(h)}_m(x_1)\|^2 dx_1   \leq \frac{C}{h^2}\int_{\Omega^h} \dist^2 (\nabla \yb^{(h)},\SO(3)) dx,
\end{equation}
where $\Rbb_m^{(h)}:(-L_1,L_2) \to \SO(3)$ is defined by
\begin{equation}
\Rbb^{(h)}_m (x_1)=\left\{\begin{array}{ll}\Rbb^{(h)}(- x_1\tb_1),& \textrm{if} \ x_1 \in (-L_1,0]\\ \Rbb^{(h)}(x_1 \tb_2),& \textrm{if} \ x_1 \in (0,L_2)  \end{array} \right. .
\end{equation}
By iterative application of (\ref{uzasno1}) and using the inequality $(x_1+\cdots+x_n)^2 \leq n(x_1^2+\cdots+x_n^2)$ and the assumption (\ref{ogranicenost})
for every open interval $I'$ compactly contained in $(-L_1,L_2)$ and $\xi \in \ZR$,  which satisfies $|\xi| \leq \dist(I',\{-L_1,L_2 \})$, we have
\begin{equation}\label{uzasno}
 \int_{I'} \|\Rbb^{(h)}_m(x_1+\xi)-\Rbb^{(h)}_m(x_1)\|^2 dx_1   \leq C_4\left(\frac{|\xi|}{h}+1 \right)^2 \frac{1}{h^2}\int_{\Omega^h} \dist^2 (\nabla \yb^{(h)},\SO(3)) dx \leq C_5(|\xi|+h)^2,
\end{equation}
Note here that the factor $(\frac{|\xi|}{h}+1)^2$ is the upper estimate of the number of terms by which the left hand side of (\ref{uzasno}) has to be estimated.
Using the Fr\'echet-Kolmogorov (see \cite[Theorem 2.21, Theorem 2.22]{Adams})  criterion, one can deduce from this that for any sequence $h_j \to 0$ there exist a subsequence $(\Rbb^{(h_{j_{1,2}})}_m)$ strongly converging in $L^2(-L_1,L_2)$ to some $\oRbb \in L^2(-L_1,L_2)$ with $\oRbb(x_1) \in \SO(3)$ for a.e. $x_1 \in (-L_1,L_2)$.
We define
$\oRbb_1:(0,L_1) \to \SO(3), \oRbb_2:(0,L_2) \to \SO(3)$ as
\begin{eqnarray*}
&&\oRbb_1(x_1)=\oRbb(-x_1),\qquad \textrm{if} \ x_1 \in (-L_1,0), \\
&&\oRbb_2(x_1)=\oRbb(x_1),\qquad \textrm{if} \ x_1 \in (0,L_2).
\end{eqnarray*}

We shall prove that $\oRbb \in W^{1,2}((-L_1,L_2);\ZR^{3 \times 3})$. Using the estimate (\ref{uzasno}) and letting $h \to 0$ we obtain that
for every $I'$ compactly contained in $(-L_1,L_2)$ and every $\xi$ which satisfies $|\xi| \leq \dist(I',\{-L_1,L_2 \})$
there exists a constant $C$ independent of $I'$ and $\xi$ such that
\begin{equation} \int_{I'} \frac{\|\oRbb(x_1+\xi)-\oRbb(x_1)\|^2}{|\xi|^2} dx_1  \leq C.
\end{equation}
>From standard theorems we obtain that $\oRbb \in W^{1,2}((-L_1,L_2); \ZR^{3 \times 3})$. This is equivalent to the fact that $\oRbb_1 \in W^{1,2}((0,L_1); \ZR^{3 \times 3})$, $\oRbb_2 \in W^{1,2}((0,L_2); \ZR^{3 \times 3})$ and $\oRbb_1(0)=\oRbb_2(0)$ in the sense of traces. In the same way one can take cylinders $C_1^h$ and $C_i^h$ for $i=3,\dots, n$ (by choosing every time a subsequence $\oRbb^{(h_{j_{1,\dots,i}})}$ of already chosen sequence  $\oRbb^{(h_{j_{1,\dots,i-1}})}$)  we have the existence of $\oRbb_i$. Moreover, the definition of $\oRbb_1$ is not ambiguous and $\oRbb_1(0)=\oRbb_2(0)=\cdots=\oRbb_n(0)$.
$$
\frac{1}{h^2}\sum_{i=1}^n \int_{C_i^h} \| \nabla \yb^{(h)}(x) - \oRbb_i(x\cdot \tb_i) \|^2 dx \leq \frac{2}{h^2} \int_{\Omega^h} \|\nabla \yb^{(h)} -\Rbb^{(h)}\|^2dx + \frac{2}{h^2} \sum_{i=1}^n \int_{C^h_i} \|\oRbb_i(x\cdot \tb_i) -\Rbb^{(h)}(x)\|^2dx
$$
Using the estimate (\ref{jjs}) and $\Rbb^{(h)} \to \oRbb_i$ in $L^2(0,L_i)$ we  obtain that
$$
\lim_{h \to 0} \frac{1}{h^2}\sum_{i=1}^n \int_{C_i^h} \| \nabla \yb^{(h)}(x) - \oRbb_i(x\cdot \tb_i) \|^2 dx = 0.
$$
\end{prooof}

\section{$\Gamma$-convergence}
\setcounter{equation}{0}

In the Theorem~\ref{glavni} we obtained the asymptotics of $\nabla \yb^{(h)}$. Still, as we are in the pure traction case, in order to obtain the asymptotics of $\yb^{(h)}$ one needs to control the constant. Thus we additionally assume that the mean value at the end of the first rod behaves nicely. Then
we obtain that in the limit in the junction point displacements from different rods have to be equal.
\begin{lemma} \label{lemasranje}
Let $(h^j)$ be a sequence that converges to $0$  and $(\yb^{(h_j)}) \subset W^{1,2}(\Omega^{h_j};\ZR^3)$ such that
\begin{equation}\label{la1}
\limsup_{j\to \infty}\frac{1}{h_j^2} \int_{\Omega^{h_j}} \|\nabla \yb^{(h_j)}\|^2 dx < \infty.
\end{equation}
Let there exist $\yb_i^0 \in W^{1,2} ((0,L_i);\ZR^3)$ such that $\yb_i^0(0)=0$ in the sense of traces and let us suppose that for every $i$,
\begin{equation}\label{la2}
\lim_{j \to \infty} \int_{(h_j,L_i) \times S_i} \| \nabla
\Reshji (\yb^{(h_j)} \circ \Qbb_i \circ\Pbhj) - (\ \EEij (\yb_{i}^0)'\ | \ 0 \ | \ 0\ )\|^2
dx =0.
\end{equation}
Let us also suppose that there exists
\begin{equation}\label{la3}
\lim_{j\to \infty} \fint_{\{L_1\} \times S_1} \yb^{(h_j)} \circ \Qbb_1 \circ \Pbhj \, dx:=C_{L_1}  \in \ZR^3.
\end{equation}
Then for $C_0:=C_{L_1}-\yb_1^0 (L_1)$  we have
\begin{equation} \label{tragovi}
\lim_{j \to \infty} \| \Reshjj\yb^{(h_j)}\circ \Qbb_1\circ\Pbhj -C_0\|_{L^2({\{h_j\} \times S_1})} = \dots=\lim_{j \to \infty} \| \Reshjn \yb^{(h_j)}\circ \Qbb_n \circ\Pbhj -C_0\|_{L^2({\{h_j\} \times S_n})}=0
\end{equation}
and
$$
\lim_{j\to\infty}\sum_{i=1}^n \|\Reshji  \yb^{(h_j)} \circ \Qbb_i \circ\Pbhj -\EEij \yb^c_i\|_{W^{1,2} ((h_j,L_i)\times S_i)}=0,
$$
where $\yb^c_i(x_1)=C_0+\yb_i^0(x_1)$.
\end{lemma}
\begin{prooof}
By applying the Poincare inequality (see part b) of \cite[Theorem 6.1-8]{Ciarlet}) to the cylinders $(0,1) \times S_i$ we have that there exists a constant $K_1$ such that for every $i\in\{1,\ldots,n\}$ and every $\yb \in W^{1,2}((0,1) \times S_i;\ZR^3)$ one has
$$
\left\|\yb-\fint_{\{1\} \times S_i} \!\!\yb\, dx\right\|_{L^{2}((0,1) \times S_i;\ZR^3)}  \leq K_1 \|\nabla \yb \|_{L^2((0,1) \times S_i;\ZR^3)}.
$$
By applying this estimate on functions of the form $\wtyb (x)=\yb ((L_i-h)x_1+h,x_2,x_3)$ we obtain that there is a constant $K_2 = \max\{ 1,L_i-h \}K_1$ such that for all $i \in \{1,\ldots, n\}$ and for all $h>0$ (small enough) and all $\yb \in W^{1,2}((h,L_i) \times S_i;\ZR^3)$ one has
\begin{equation} \label{druga}
\left\|\yb-\fint_{\{L_i\} \times S_i} \!\!\yb\, dx\right\|_{L^{2}((h,L_i) \times S_i;\ZR^3)}
\leq K_2 \|\nabla \yb \|_{L^2((h,L_i) \times S_i;\ZR^3)}.
\end{equation}
In a similar way we obtain
\begin{equation} \label{drugaipo}
\left\|\yb-\fint_{\{h\} \times S_i} \!\!\yb\, dx\right\|_{L^{2}((h,L_i) \times S_i;\ZR^3)}  \leq K_2' \|\nabla \yb \|_{L^2((h,L_i) \times S_i;\ZR^3)}.
\end{equation}
Moreover, by using the same rescaling of the domain, from continuity of traces we obtain that there is a constant $K_3$ such that for all $i,h$ and $\yb \in W^{1,2}((h,L_i) \times S_i;\ZR^3)$ one has
\begin{equation} \label{tragovi1}
\|\yb\|_{L^2(\{h\} \times S_i;\ZR^3)}+\|\yb\|_{L^2(\{L_i\} \times S_i;\ZR^3)} \leq K_3 \| \yb\|_{W^{1,2}((h,L_i)\times S_i ;\ZR^3)}.
\end{equation}

By applying the Poincare inequality (of the same form as before) to the domain $T$ on functions given by $\wtyb (x)=\yb(hx)$ we have that there exist a constant $K_4$ such that for all $i,h$ and $\yb \in W^{1,2}(T^h;\ZR^3)$ one has
$$
\left\|\yb-\fint_{\Qbb_i(\{h\}\times hS_i)} \!\!\yb\, dx\right\|_{L^{2}(T^h;\ZR^3)} \leq hK_4 \| \nabla \yb\|_{L^{2}(T^h;\ZR^3)}.
$$
In the similar way as before we conclude that there exists a constant $K_5 = 2 K_4$ such that for all $i,l,h$ and $\yb \in W^{1,2}(T^h;\ZR^3)$ one has
\begin{equation} \label{peta}
\left\|\fint_{\Qbb_i(\{h\} \times hS_i)} \!\!\yb\, dx-\fint_{\Qbb_l(\{h\} \times hS_l)} \!\!\yb\, dx\right\| \leq  \frac{K_5}{\sqrt{h} } \| \nabla \yb \|_{L^{2}(T^h;\ZR^3)}. \end{equation}

We now apply inequality (\ref{druga}) to the sequence $\Reshjj \yb^{(h_j)} \circ \Qbb_1 \circ\Pbhj -\EEjj\yb^c_1$ to obtain
\begin{eqnarray*}
&&\|\Reshjj  \yb^{(h_j)} \circ \Qbb_1 \circ\Pbhj -\EEjj \yb^c_1- (\fint_{\{L_1\} \times S_1} \yb^{(h_j)}\circ \Qbb_1 \circ \Pbhj \, dx-C_{L_1})\|_{L^{2} ((h_j,L_1)\times S_1;\ZR^3)}\\ && \hspace{10ex} \leq
K_2 \| \nabla \Reshjj (\yb^{(h_j)} \circ \Qbb_1 \circ\Pbhj) - (\ \EEjj ( \yb_{1}^0)' \ | \ 0 \ | \ 0\ )\|_{L^2((h_j,L_1) \times S_1; \ZR^3)}.
\end{eqnarray*}
Now, using the assumptions (\ref{la2}) and (\ref{la3}) we obtain that $\|\Reshjj  \yb^{(h_j)} \circ \Qbb_1 \circ\Pbhj -\EEjj \yb^c_1\|_{W^{1,2} ((h_j,L_1)\times S_1;\ZR^3)} \to 0$. The estimate (\ref{tragovi1}) now  implies
\begin{equation}\label{sesta}
\lim_{j \to \infty}\| \Reshjj \yb^{(h_j)}\circ \Qbb_1 \circ\Pbhj -C_0\|_{L^2(\{h_j\}\times S_1;\ZR^3)}=0.
\end{equation}
By applying (\ref{peta}) for $l=1$ and $i \neq 1$ to the sequence $\yb^{(h_j)}$ we obtain
$$
\left\|\fint_{\Qbb_i(\{h_j\} \times h_jS_i)} \!\!\yb^{(h_j)}\, dx-\fint_{\Qbb_1(\{h_j\} \times h_jS_1)} \!\!\yb^{(h_j)}\, dx\right\| \leq  \frac{K_5}{\sqrt{h_j} } \| \nabla \yb^{(h_j)} \|_{L^{2}(T^{h_j};\ZR^3)}.
$$
Now we change the variables in the integrals on the left hand side (also note that $T^{h_j} \subset \Omega^{h_j}$) to obtain
$$
\left\|\fint_{\{h_j\} \times S_i} \!\!\yb^{(h_j)} \circ \Qbb_i \circ \Pbhj\, dx-\fint_{\{h_j\} \times S_1} \!\!\yb^{(h_j)}\circ \Qbb_1 \circ \Pbhj\, dx\right\| \leq  \frac{K_5}{\sqrt{h_j} } \| \nabla \yb^{(h_j)} \|_{L^{2}(\Omega^{h_j};\ZR^3)}.
$$
Therefore (\ref{la1}) and (\ref{sesta}) imply $\fint_{\{h_j\} \times S_i} \Reshji \yb^{(h_j)}\circ \Qbb_i \circ\Pbhj \,dx \to C_0$ for all $i$'s.
By applying the inequality (\ref{drugaipo}) to the sequence
$\Reshji \yb^{(h_j)}\circ \Qbb_i \circ\Pbhj -\EEij \yb^c_i$ for $i \neq 1$ we obtain that  $\|\Reshji  \yb^{(h_j)} \circ \Qbb_i \circ\Pbhj -\EEij \yb^c_i\|_{W^{1,2} ((h_j,L_i)\times S_i;\ZR^3)} \to 0$.
Then  (\ref{tragovi}) follows immediately from (\ref{tragovi1}) for $i=1$ and using the fact that $\|\EEij \yb^c_i-C_0\|_{L^{2} (\{h_j\} \times S_i;\ZR^3)} = |S_i|^{1/2}\|\EEij \yb^c_i(h_j)-C_0\|\to 0$.
\end{prooof}

In the following we use the notation
$$
(\yb,\db^2,\db^3) = ((\yb_1,\db^2_1,\db^3_1), \dots, (\yb_n,\db^2_n,\db^3_n))
$$
to collect deformations of all rods.

Combining the results of Theorem~\ref{glavni} and Lemma~\ref{lemasranje} we obtain the following result.
\begin{corollary}\label{ct3.1l4.1}
Let $(\yb^{(h)}) \subset W^{1,2}(\Omega^h;\ZR^3)$ be such that
\begin{equation}\label{ogr}
\limsup_{h \to 0} \frac{1}{h^4} \int_{\Omega^h} \dist^2 (\nabla \yb^{(h)} ,\SO(3))dx < +\infty,
\end{equation}
\begin{equation}\label{la3a}
\lim_{j\to \infty} \fint_{\{L_1\} \times S_1} \yb^{(h_j)} \circ \Qbb_1 \circ \Pbhj \, dx:=C_{L_1}  \in \ZR^3.
\end{equation}
Then for every sequence in $\ZR_+$ converging to 0 there exist a subsequence $(h_j)$ and $\yb_i \in W^{2,2}((0,L_i); \ZR^3)$, $\db^2_i, \db^3 \in W^{1,2}((0,L_i); \ZR^3)$ such that for $\Rbb_i = (\yb_i' \ \db^2_i \ \db^3_i)$ one has
\begin{eqnarray*}
&&(\yb, \db^2, \db^3) \in \mathcal{A}:=\{ ((\yb_1,\db^2_1,\db^3_1),\dots,(\yb_n,\db^2_n,\db^3_n) ) \in (W^{2,2}((0,L_1);\ZR^3) \times W^{1,2}((0,L_1);\ZR^3) \times W^{1,2}((0,L_1);\ZR^3))\\ && \hspace{20ex} \times \dots \times (W^{2,2}((0,L_n);\ZR^3) \times W^{1,2}((0,L_n);\ZR^3) \times W^{1,2}((0,L_n);\ZR^3)): \\
 && \hspace{20ex} \Rbb_i
 \in \SO(3) \textrm{ a.e. and }  \yb_1(0)=\cdots=\yb_n(0), \Rbb_1(0)\Qbb_1^T=\cdots=\Rbb_n(0)\Qbb_n^T
 \}
\end{eqnarray*}
and
\begin{equation}\label{kvgC}
\lim_{j \to \infty} \frac{1}{h_j^2} \sum_{i=1}^n \|\yb^{(h_j)}\circ \Qbb_i- D_i(\yb_i,\db^2_i,\db^3_i)\|^2_{W^{1,2}((h_j,L_i) \times h_jS_i;\ZR^3)}=0,
\end{equation}
where $D_i(\yb_i,\db^2_i,\db^3_i)(x_1,x_2,x_3)=\yb_i(x_1)+x_2 \db^2_i(x_1)+x_3 \db^3_i (x_1)$, for $x \in (h_j,L_i) \times h_j S_i$.
\end{corollary}
\begin{prooof}
>From (\ref{ogr}) it follows that the assumption of Theorem~\ref{glavni} is fulfilled. Therefore there exist a subsequence $(h_j)$ converging to 0 and
$\oRbb_i \in W^{1,2}((0,L_i), \SO(3))$, $i=1,\dots, n$ such that $\oRbb_1(0)=\oRbb_2(0)=\cdots=\oRbb_n (0)$ in the sense of traces and
$$
\lim_{j \to \infty} \frac{1}{h_j^2} \sum_{i=1}^n\int_{C_i^{h_j}} \| \nabla \yb^{(h_j)}(x) - \oRbb_i(x \cdot \tb_i) \|^2 dx = 0.
$$
We rewrite this convergence to obtain
\begin{eqnarray}
\nonumber
0&=& \lim_{j \to \infty} \frac{1}{h_j^2} \sum_{i=1}^n\int_{(h_j,L_i)\times h_j S_i} \| \nabla \yb^{(h_j)}(\Qbb_i x) - \oRbb_i(\Qbb_i x \cdot \tb_i) \|^2 dx\\
\nonumber
&=& \lim_{j \to \infty} \frac{1}{h_j^2} \sum_{i=1}^n\int_{(h_j,L_i)\times h_j S_i} \| \nabla (\yb^{(h_j)}\circ \Qbb_i )(x) \Qbb_i^T - \oRbb_i(x \cdot \eb_1) \|^2 dx\\
&=& \lim_{j \to \infty} \frac{1}{h_j^2} \sum_{i=1}^n\int_{(h_j,L_i)\times h_j S_i} \| \nabla (\yb^{(h_j)}\circ \Qbb_i )(x) - \oRbb_i(x_1)\Qbb_i \|^2 dx. \label{neka}
\end{eqnarray}
Now we define
\begin{eqnarray*}
&& (\yb^0_i)' = \oRbb_i(x_1)\Qbb_i \eb_1, \qquad \yb^0_i(0)=0,\\
&& \db^2_i = \oRbb_i(x_1)\Qbb_i \eb_2, \\
&& \db^3_i = \oRbb_i(x_1)\Qbb_i \eb_3, \\
&& \Rbb_i = \oRbb_i(x_1)\Qbb_i.
\end{eqnarray*}
Since $\oRbb_i \in W^{1,2}((0,L_i), \SO(3))$ it follows that $\yb^0_i \in W^{2,2} ((0,L_i); \ZR^3)$ and $\Rbb_i \in W^{1,2}((0,L_i), \SO(3))$. By the trace property of $\oRbb_i$ we obtain $\Rbb_1(0)\Qbb_1^T=\cdots=\Rbb_n(0)\Qbb_n^T$.

In the sequel we want to apply the Lemma~\ref{lemasranje}. Therefore we check the assumptions of the lemma. First, we estimate the norm of a matrix by  the distance of the matrix to $\SO(3)$  and the norm of an arbitrary rotation to obtain
$$
\int_{\Omega^{h_j}} \|\nabla \yb^{(h_j)}\|^2 dx \leq 2\int_{\Omega^{h_j}}
\dist^2 (\nabla \yb^{(h_j)}(x), \SO(3))dx+ 
C h_j^2.
$$
Using (\ref{ogr}) we obtain that (\ref{la1}) is fulfilled.

Changing the coordinates in (\ref{neka}) we obtain
\begin{eqnarray*}
\nonumber
0
&=& \lim_{j \to \infty}  \sum_{i=1}^n\int_{(h_j,L_i)\times S_i} \| \nabla (\yb^{(h_j)}\circ \Qbb_i \circ \Pbhj)(x) \nabla \Pbhj  - \oRbb_i(x_1)\Qbb_i \|^2 dx.
\end{eqnarray*}
This implies that (\ref{la2}) is fulfilled with $\yb^0_i$ defined above. The assumption (\ref{la3}) is fulfilled by (\ref{la3a}). Therefore we can apply Lemma~\ref{lemasranje} to obtain that for $C_0:=C_{L_1}-\yb_1^0 (L_1)$  we have (\ref{tragovi})
and
\begin{equation}\label{l2}
\lim_{j\to\infty}\sum_{i=1}^n \|\Reshji  \yb^{(h_j)} \circ \Qbb_i \circ\Pbhj -\EEij \yb_i\|_{W^{1,2} ((h_j,L_i)\times S_i)}=0,
\end{equation}
where $\yb_i(x_1)=C_0+\yb_i^0(x_1)$. Since $\yb_i' = (\yb^0_i)'$ from (\ref{l2}) and (\ref{neka}) we obtain (\ref{kvgC}).
>From (\ref{kvgC}) and the estimate
$$
\| \Reshji \yb^{(h_j)}\circ \Qbb_i \circ\Pbhj -\yb_i\|_{L^2({\{h_j\} \times S_i})} \leq C \| \Reshji\yb^{(h_j)}\circ \Qbb_i \circ\Pbhj -\yb_i\|_{W^{1,2}({(h_j,L_i) \times S_i})}
$$
(for details see the proof of Lemma~\ref{lemasranje}) we obtain
\begin{equation}\label{kvg2}
\lim_{j \to \infty} \| \Reshji \yb^{(h_j)}\circ \Qbb_i \circ\Pbhj -\yb_i\|_{L^2({\{h_j\} \times S_i})} = 0, \qquad i=1\ldots, n.
\end{equation}
Now, (\ref{tragovi}) and (\ref{kvg2}) imply
$$
|S_i|^{1/2}\|\yb_i(h_j)-C_0\|=\|\yb_i-C_0\|_{L^2({\{h_j\} \times S_i})} \to 0,
$$
for all $i=1, \ldots, n$. This implies that $\yb_1(0)=\dots=\yb_n(0)=C_0$.

Thus we obtain that $(\yb, \db^2,\db^3) \in \cA$.
\end{prooof}

\begin{remark}\label{rDi}
The structure of the functions $D_i(\yb_i, \db^2_i, \db^3_i)$ defined after (\ref{kvgC}) is essentially  one--dimensional. It stands as the limit displacement for the $i$-th rod. The function $\yb_i$ describes deformation of the middle curve of the $i$-th rod, while the vectors $\db^2_i$ and $\db^3_i$ span the normal plane of the deformed middle curve (since $\Rbb_i = (\yb_i' \ \db^2_i \ \db3_i) \in \SO(3)$). Since the rod is assumed thin, variables $x_2$ and $x_3$ (cross-sectional coordinates of $h S_i$) are of order $h$ so the terms involving these terms can be considered as first correctors to the leading order approximation $\yb_i$ of the $i$-th rod. Note as well that the convergence (\ref{kvgC}) will be the one which will be used to formulate the asymptotics of the infimizing sequence.
\end{remark}

\begin{proposition}\label{propo}
Let the functional $I$ is defined by
$$
I(\yb,\db^2,\db^3)=\left\{ \begin{array}{ll} \displaystyle \sum_{i=1}^n \frac{1}{2} \int_0^{L_i} \Qu^i_2 (\Rbb_i^T \Rbb_{i}') dx_1 & \textrm{if} \ (\yb,\db^2,\db^3) \in \mathcal{A}, \\   +\infty & \textrm{otherwise}   \end{array}\right. $$
    where $\Rbb_i:=(\yb_{i}', \db^2_i,\db^3_i)$, while the class $\mathcal{A}$ is given by
\begin{eqnarray*}
 && \mathcal{A}:=\{ ((\yb_1,\db^2_1,\db^3_1),\dots,(\yb_n,\db^2_n,\db^3_n) ) \in (W^{2,2}((0,L_1);\ZR^3) \times W^{1,2}((0,L_1);\ZR^3) \times W^{1,2}((0,L_1);\ZR^3))\\ && \hspace{20ex} \times \dots \times (W^{2,2}((0,L_n);\ZR^3) \times W^{1,2}((0,L_n);\ZR^3) \times W^{1,2}((0,L_n);\ZR^3)): \\
 && \hspace{20ex} \Rbb_i
 \in \SO(3) \textrm{ a.e. and }  \yb_1(0)=\cdots=\yb_n(0), \Rbb_1(0)\Qbb_1^T=\cdots=\Rbb_n(0)\Qbb_n^T
 \}.
\end{eqnarray*}
The quadratic forms $\Qu^i_2:\mathbb{M}^{3 \times 3}_{skew} \to [0,\infty)$ are defined by
 \begin{equation} \label{minimum}
 \Qu^i_2(\Abb):=\min_{\alpha \in W^{1,2}(S_i;\ZR^3)} \int_{S_i} \Qu^i_3 \left( \Abb \left( \begin{array}{c} 0 \\ x_2 \\ x_3 \end{array}\right) \Bigg| \partial_2 \alpha_{} \Bigg| \partial_3 \alpha_{} \right) dx_2 dx_3,
 \end{equation}
where
$$
\Qu^i_3 (\Gbb)= \frac{\partial^2 W}{\partial \Fbb^2} (\Qbb_i^T)(\Gbb\Qbb_i^T,\Gbb\Qbb_i^T).
$$
Then the following two statements hold.
\begin{itemize}
\item(liminf inequality)
Let $\yb_i \in  W^{1,2} ((0,L_i);\ZR^3)$, $\db_i^2,\db_i^3 \in L^2((0,L_i);\ZR^3)$. Then for every sequence $(h_j)\subset (0,\infty)$ converging to $0$ and every sequence $(\yb^{(h_j)}) \subset W^{1,2}(\Omega^{h_j};\ZR^3)$ such that
    \begin{equation}\label{kvg}
    \lim_{j \to \infty} \frac{1}{h_j^2} \sum_{i=1}^n \|\yb^{(h_j)}\circ \Qbb_i- D_i(\yb_i,\db^2_i,\db^3_i)\|^2_{W^{1,2}((h_j,L_i) \times h_jS_i;\ZR^3)}=0,
    \end{equation}
    where $D_i$ are defined in (\ref{kvgC}), 
we have that
$$
I(\yb,\db^2,\db^3) \leq \liminf_{j \to \infty} \frac{1}{h_j^4} E^{(h_j)} (\yb^{(h_j)}); $$
\item (limsup inequality) For every sequence  $(h_j)\subset (0,\infty)$ converging to $0$ and for every $\yb_i \in W^{1,2} ((0,L_i);\ZR^3)$, $\db^2_i,\db^3_i \in L^2((0,L_i);\ZR^3)$  there exist a sequence $(\yb^{(h_j)}) \subset W^{1,2} (\Omega^{h_j};\ZR^3)$ such that
    $$
    \lim_{j \to \infty} \frac{1}{h_j^2} \sum_{i=1}^n \|\yb^{(h_j)}\circ \Qbb_i- D_i(\yb_i,\db^2_i,\db^3_i)\|^2_{W^{1,2}((h_j,L_i) \times h_jS_i;\ZR^3)}=0
    $$
    and
    $$
    \lim_{j \to \infty} \frac{1}{h_j^4} E^{(h_j)}(\yb^{(h_j)})=I(\yb,\db^2,\db^3).
    $$
\end{itemize}
\end{proposition}

\begin{remark} \label{napp}
As it is noted in Remark 3.4. in \cite{Muller1} each minimization problem in (\ref{minimum}) has a solution and it can be equivalently computed
on the class of functions
$$
V_i= \{ \alb \in W^{1,2} (S_i;\ZR^3): \int_{S_i} \alb dx_2 dx_3=\int \nabla \alb dx_2 dx_3 =0 \}.
$$
It can be also shown that for every $i$  the minimizer is unique in $V_i$ and that minimizer in $V_i$ depends linearly on the entries $(a_{ij})$ of $A$. Hence $\Qu^i_2$ is in fact a quadratic form of $A$. In the isotropic case ($W(\Fbb)=W(\Fbb \Rbb)$ for every $\Fbb \in \mathbb{M}^{3 \times 3}$ and $\Rbb \in \SO(3)$) for every $i$ we have $\Qu^i_3 (\Gbb)= \frac{\partial^2 W}{\partial \Fbb^2} (\Id)(\Gbb,\Gbb)$. In this case there are also some explicit formulas for $\Qu^i_2$ (see Remarks 3.5. and 3.6. in \cite{Muller1}).
\end{remark}
\begin{prooof}
Let us first prove the $\liminf$ inequality.

Let $(\yb,\db^2,\db^3) \in \mathcal{A}$ and let
$0<h_j\to 0$ and $\yb^{(h_j)} \subset W^{1,2}(\Omega^{h_j}; \ZR^3)$ satisfy (\ref{kvg}).
Let us also fix $\delta >0$.  Then, after rescaling each convergence in the sum (\ref{kvg}) to the fixed domain $(\delta, L_i) \times S_i$ we obtain that for every $i\in \{1, \ldots, n\}$ one has
\begin{eqnarray*}
&&\| \Reshji \yb^{(h_j)} \circ \Qbb_i \circ \Pbhj- \EEij \yb_i\|_{W^{1,2} ((\delta,L_i) \times S_i;\ZR^3)} \to 0, \\
&&\|(\frac{1}{h_j}\partial_2 (\Reshji\yb^{(h_j)}\circ \Qbb_i\circ\Pbhj),\frac{1}{h_j} \partial_3(\Reshji \yb^{(h_j)}\circ \Qbb_i\circ\Pbhj)) - (\EEij \db_i^2, \EEij \db_i^3)\|_{L^2((\delta,L_i)\times S_i;\ZR^3 \times \ZR^3)} \to 0.
\end{eqnarray*}
Now, by using Theorem 3.1. from \cite{Muller1} on each rod separately (applying it to the energy density functions $W^{\Qbb_i^T}(\Fbb):=W(\Fbb\Qbb_i^T)$) we conclude that for every $\delta$ and for every $i$ we have
  \begin{eqnarray*} \frac{1}{2} \int_{\delta}^{L_i} \Qu^i_2 (\Rbb_i^T \Rbb_{i}') dx_1 &\leq& \liminf_{j\to \infty} \frac{1}{h_j^4} \int_{(\delta, L_i)\times S_i} W^{\Qb_i^T} (\nabla_{h_j} \Reshji (\yb^{(h_j)} \circ \Qbb_i \circ\Pbhj)) dx\\
  &=&\liminf_{j \to \infty} \frac{1}{h_j^4} \int_{(\delta,L_i) \times h_j S_i} W(\nabla \yb_i^{(h_j)}(\Qbb_i x)\Qbb_i \Qbb_i^T)dx \\
  &=&\liminf_{j \to \infty} \frac{1}{h_j^4} \int_{\Qbb_i((\delta,L_i) \times h_j S_i)} W(\nabla \yb_i^{(h_j)} (x)) dx, \end{eqnarray*}
where we have used the notation $\nabla_h = (\partial_1 \ \frac{1}{h}\partial_2 \ \frac{1}{h}\partial_3)$.
By summing all these inequalities we obtain that for every $\delta>0$ one has
$$
\sum_{i=1}^n \frac{1}{2} \int_{\delta}^{L_i} \Qu^i_2 (\Rbb_i^T \Rbb_{i}') dx_1 \leq \liminf_{j \to \infty} \frac{1}{h_j^4} E^{(h_j)} (\yb^{(h_j)}).
$$
By letting $\delta \to 0$ we obtain
$$
I(\yb,\db^2,\db^3) \leq \liminf_{j \to \infty} \frac{1}{h_j^4} E^{(h_j)} (\yb^{(h_j)}).
$$

Let us now suppose that $(\yb,\db^2,\db^3) \notin \mathcal{A}$. We
have to see that for every sequence $(\yb^{(h_j)})\subset
W^{1,2}(\Omega^{h_j};\ZR^3)$ such that (\ref{kvg}) holds one has
$\liminf_{j \to \infty} \frac{1}{h_j^4} E^{(h_j)} (\yb^{(h_j)})$
$=+ \infty$. Let us suppose the opposite, i.e., $\liminf_{j \to
\infty} \frac{1}{h_j^4} E^{(h_j)} (\yb^{(h_j)})<  \infty$. Using
the property of the stored energy function $W$ we estimate
\begin{equation}\label{bded}
C_W \frac{1}{h_j^4} \int_{\Omega^{h_j}} \dist^2(\nabla
\yb^{(h_j)}(x),\SO(3))dx \leq \frac{1}{h_j^4} E^{(h_j)}
(\yb^{(h_j)}) < \infty.
\end{equation}
>From the convergence (\ref{kvg}) one can easily conclude, using the continuity of the trace operator and the fact we can control the change of the domain (similarly as in Lemma~\ref{lemasranje} and Corollary~\ref{ct3.1l4.1}), that
$$ \lim_{j\to \infty} \fint_{\{L_1\} \times S_1} \yb^{(h_j)} \circ \Qbb_1 \circ \Pbhj \, dx:=\yb_1(L_1). $$
Thus the assumptions of Corollary~\ref{ct3.1l4.1} are fulfilled and we can conclude, by the uniqueness of the limit, that $(\yb,\db^2,\db^3) \in \mathcal{A}$, which is a contradiction.

To prove $\limsup$ inequality we have to construct the appropriate
sequence. Let us take $(\yb,\db^2,\db^3) \in \mathcal{A}$. Let us
in addition suppose $\yb_i \in C^2 ([0,L_i];\ZR^3)$, $\db_i^2,
\db_i^3 \in C^1([0,L_i];\ZR^3)$ (note that $\yb_i \in C^2
([0,L_i];\ZR^3)$ is an immediate consequence of $\db_i^2, \db_i^3
\in C^1([0,L_i];\ZR^3)$ and $\Rbb_i \in \SO(3)$). Let us define
$\yb^{(h_j)}$ in the following way
\begin{eqnarray*}
\yb^{(h_j)}(x) &=& \yb_i(0)+ \Rbb_i(0) \Qbb_i^T x, \quad
\textrm{for } x\in T^{h_j}\ \textrm{(the definition is not ambiguous)}, \\
\yb^{(h_j)}(\Qbb_i \Pbhj(x)) &=&
\yb_i(x_1-h_j)+h_j\yb_{i}'(0)+ h_j x_2 ( \db_i^2
(x_1-h_j)-\alpha^{(h_j)}(x_1)(\db_{i}^2)'(0))\\ & &+ h_j x_3
(\db_i^3(x_1-h_j)-\alpha^{(h_j)}(x_1)(\db_{i}^3)'(0))+h_j^2 \betab_i^j
(x),\quad  \textrm{for} \ x \in
(h_j,L_i) \times S_i.
\end{eqnarray*}
where
$ \alpha^{(h_j)} \in C^1([h_j,L_i]; \ZR^3) $ are such that
\begin{eqnarray*}
&&\alpha^{(h_j)}(h_j)=0, \quad (\alpha^{(h_j)})'(h_j)=1, \quad (\alpha^{(h_j)})'(2h_j)=0,\\
&&\alpha^{(h_j)}(x_1)=0, \mbox{ for } x_1 \geq 2h_j, \quad \sup_{j}
\|\alpha^{(h_j)}\|_{\infty}<Ch_j,\quad \sup_{j}
\|(\alpha^{(h_j)})'\|_{\infty}<\infty
\end{eqnarray*}
(e.g.
$\alpha^{(h_j)}(x_1)=\frac{1}{h_j^2}
x_1^3-\frac{5}{h_j}x_1^2+8x_1-4h_j$, for $x_1 \in [h_j,2h_j]$ and
$0$ otherwise).
The functions $\betab_i^j:[0,L_i] \times S_i \to \ZR^3$ are chosen
such that $\betab_i^j(x)=\gamma_i^j(x_1)\betab_i(x)$,
where
$\gamma_i^j \in C^1([0,L_i];\ZR)$ are such that
\begin{eqnarray*}
&&\gamma_i^j(x_1)=0, \mbox{ for } x_1 \leq h_j, \quad \gamma_i^j(x_1)=1, \mbox{ for } x_1 \geq 2h_j, \quad \|
\gamma_i^j \|_{\infty}<C, \quad \| (\gamma_i^j)' \|_{\infty}<
\frac{C}{h_j}
\end{eqnarray*}
(e.g. $\gamma_i^j(x) =
-\frac{2}{h_j^3}x^3+\frac{9}{h_j^2}x^2-\frac{12}{h_j}x+5$) and
$\betab_i \in C^1([0,L_i] \times S_i; \ZR^3)$.
Then we have
\begin{eqnarray*}
\nabla \yb^{(h_j)}(x) &=&  \Rbb_i(0) \Qbb_i^T, \qquad \textrm{for } x \in T^{h_j}, \\
\nabla \yb^{(h_j)}(\Qbb_i \Pbhj(x))\Qbb_i &=&
\Rbb_i(x_1-h_j) -(0 \ | \alpha^{(h_j)}(x_1)(\db_{i}^2)'(0) \ | \
{\alpha^{(h_j)} (x_1)}{} (\db_{i}^3)'(0) )1_{h_j \leq x_1 \leq 2h_j}\\
& &+ h_j (x_2 (\db_{i}^2)'(x_1-h_j)+x_3 (\db_{i}^3)' (x_1-h_j)\ |\
\partial_2 \betab_{i}^j(x)  \ |\
\partial_3 \betab_{i}^j(x))\\
& & +h_j (-x_2(\alpha^{(h_j)})'(x_1)(\db_{i}^2)'(0)-x_3 (\alpha^{(h_j)})'(x_1)
(\db_{i}^3)'(0)\ |\  0  \ |\ 0)1_{h_j \leq x_1 \leq 2h_j}\\
&&+h_j^2(
\partial_1 \betab_{i}^j(x) \ | \ 0 \ | \
0), \qquad \textrm{for } x\in (h_j,L_i) \times S_i.
\end{eqnarray*}
Note that $\yb^{(h_j)} \in C^1(\Omega^{h_j}; \ZR^3)$. It can be
easily seen, by the dominated convergence theorem, that for every
$i$ we have $\lim_{j\to \infty} \frac{1}{h_j}\| \yb^{(h_j)}
\circ \Qbb_i-\yb_i\|_{L^2((h_j,L_i) \times
h_j S_i;\ZR^3)}=0$ and $\lim_{j\to \infty}\frac{1}{h_j} \| \nabla
\yb^{(h_j)}\circ\Qbb_i  \Qbb_i-\Rbb_i\|_{L^2((h_j,L_i) \times h_j S_i;\ZR^3)}=0$ which together implies that $\yb^{(h_j)}$ satisfies (\ref{kvg}).
Now we have to prove the $\limsup$ inequality for this sequence.

Let us for $x \in (h_j,L_i) \times S_i$ denote
\begin{equation}
\Bbbb_i^{(h_j)}(x)=\frac{\Rbb_i(x_1-h_j)^T \nabla \yb^{(h_j)}(\Qbb_i \Pbhj(x))-\Qbb_i^T}{h_j}.
\end{equation}
Then
\begin{eqnarray*}
\Bbbb_i^{(h_j)}(x)\Qbb_i&=&
\Rbb_i(x_1-h_j)^T (x_2 (\db_{i}^2)'(x_1-h_j)+x_3 (\db_{i}^3)' (x_1-h_j)\ |\  \partial_2 \betab_{i}^j(x)  \ |\ \partial_3 \betab_{i}^j(x))\\
&&-\Rbb_i(x_1-h_j)^T(0 \ | \frac{\alpha^{(h_j)}(x_1)}{h_j} (\db_{i}^2)'(0) \ | \ \frac{\alpha^{(h_j)}(x_1)}{h_j} (\db_{i}^3)'(0) )1_{h_j \leq x_1 \leq  2h_j} \\
& &+\Rbb_i(x_1-h_j)^T(-x_2(\alpha^{(h_j)})'(x_1) (\db_{i}^2)'(0)-x_3 (\alpha^{(h_j)})'(x_1)(\db_{i}^3)'(0)\ |\  0  \ |\ 0)1_{h_j \leq x_1 \leq  2h_j}\\ & &+h_j \Rbb_i(x_1-h_j)^T( \partial_1 \betab_{i}^j(x) \ | \ 0 \ | \ 0).
\end{eqnarray*}
Note that for every $\delta>0$ one has
$$
\Bbbb_i^{(h_j)}(x)\Qbb_i \to \Rbb_i(x_1)^T (x_2 (\db_{i}^2)'(x_1)+x_3 (\db_{i}^3)' (x_1)\ |\  \partial_2 \betab_{i}(x)  \ |\ \partial_3 \betab_{i}(x)), \qquad \textrm{a.e. } x \in (\delta, L_i).
$$
For every $i$, we look at the sequence $(f^i_j)_j$ of functions $f^i_j: (0,L_i) \times S_i \to[0,+\infty)$ defined by
\begin{eqnarray*}
&& f^i_j(x)=0,\qquad \textrm{for} \ x \in (0,h_j) \times S_i, \\
&& f^i_j(x)= \frac{1}{h_j^2} W(\nabla \yb^{(h_j)}(\Qbb_i \Pbhj(x) ))=\frac{1}{h_j^2} W(\Qbb_i^T+h_j \Bbbb_i^{(h_j)}(x)),\qquad \textrm{for} \ x \in (h_j,L_i) \times S_i,
\end{eqnarray*}
where the equality in the second line holds by the objectivity of $W$.
Since $W$ is $C^2$ in the neighborhood of $\SO(3)$ and has extreme on $\SO(3)$ and $\Bbbb_i^{(h_j)}$ is bounded, by the Taylor theorem,  for every $i$, one has
$$
f^i_j (x) \to
\frac{1}{2} \Qu_3^i(\Rbb_i^T (x_2 (\db_{i}^2)'+x_3(\db_{i}^3)'\ | \ \partial_2 \betab_{i} \ | \partial_3 \betab_{i})),\quad \textrm{a.e.}\ x \in (0,L_i) \times S_i.
$$
Since $\Bbbb^{(h_j)}_i$ is bounded the sequence  $f^i_j$
is also bounded in $L^\infty((0,L_i) \times S_i;\ZR^3)$ also by the Taylor theorem. Thus by the dominated convergence theorem we have
\begin{eqnarray*}
\lim_{j \to \infty} \frac{1}{h_j^4} \int_{C_i^{h_j}} W(\nabla \yb^{(h_j)})dx&=&\lim_{j \to \infty} \int_{(h_j,L_i)\times S_i} \frac{1}{h_j^2} W(\nabla \yb^{(h_j)}(\Qbb_i \Pbhj(x))) dx \\ &=&
\lim_{j \to \infty} \int_{(0,L_i)\times S_i} f^i_j(x) dx =
\frac{1}{2} \int_{(0,L_i)\times S_i} \Qu_3^i(\Rbb_i^T (x_2 (\db_{i}^2)'+x_3(\db_{i}^3)'\ | \ \partial_2 \betab_{i} \ | \partial_3 \betab_{i})).
\end{eqnarray*}
Also note that for the chosen sequence $\yb^{(h_j)}$, for every $j$, one has $W(\nabla \yb^{(h_j)}\mid_{T^{h_j}}(x))=0$ and thus
\begin{eqnarray*}
\lim_{j \to \infty} \frac{1}{h_j^4} \int_{\Omega^{h_j}} W(\nabla \yb^{(h_j)})dx &=& \sum_{i=1}^n \lim_{j \to \infty} \frac{1}{h_j^4} \int_{C_i^{h_j}} W(\nabla \yb^{(h_j)})dx\\
&=&\sum_{i=1}^n \frac{1}{2} \int_{(0,L_i)\times S_i} \Qu_3^i(\Rbb_i^T (x_2 (\db_{i}^2)'+x_3(\db_{i}^3)'\ | \ \partial_2 \betab_{i} \ | \partial_3 \betab_{i})).
\end{eqnarray*}
Thus for $(\yb,\db^2,\db^3) \in \mathcal{A}$.
s.t. $\yb_i \in C^2 ([0,L_i];\ZR^3)$, $\db_i^2, \db_i^3 \in C^1([0,L_i];\ZR^3)$ and arbitrary $\betab_i \in C^1 ([0,L_i] \times S_i;\ZR^3)$ we have that there exists a sequence $(\yb^{(h_j)}) \subset W^{1,2} (\Omega^{h_j};\ZR^3)$ such that
$$
\lim_{j \to \infty} \frac{1}{h_j^2} \sum_{i=1}^n \|\yb^{(h_j)}\circ \Qbb_i- D_i(\yb_i,\db^2_i,\db^3_i)\|^2_{W^{1,2}((h_j,L_i) \times h_jS_i;\ZR^3)}=0
$$
and
$$
\lim_{j \to \infty} \frac{1}{h_j^4} E^{(h_j)}(\yb^{(h_j)})=\frac{1}{2} \int_{(0,L_i)\times S_i} \Qu_3^i(\Rbb_i^T (x_2 (\db_{i}^2)'+x_3(\db_{i}^3)'\ | \ \partial_2 \betab_{i} \ | \ \partial_3 \betab_{i})).
$$
Now let us now consider the general case and take an arbitrary $(\yb,\db^2,\db^3) \in \mathcal{A}$. For every $i$ we choose a sequence $(\wtRbb_i^{(j)}) \subset C^1([0,L_i]; \mathbb{M}^{3 \times 3})$ such that $\wtRbb_i^{(j)} \to  (\yb_{i}',\db_i^2,\db_i^3) = \Rbb_i$ in $W^{1,2}((0,L_i); \mathbb{M}^{3 \times 3})$. By making a slight correction, namely taking $\whRbb_i^{(j)}=\Rbb_i(0) \wtRbb_i^{(j)}(0)^{-1} \wtRbb_i^{(j)}$ (this can be done for $j$ large enough due to Sobolev embedding theorem) we also have $(\whRbb_i^{(j)}) \subset C^1([0,L_i]; \mathbb{M}^{3 \times 3})$, $\whRbb_i^{(j)} \to \Rbb_i =(\yb_{i}',\db_i^2,\db_i^3)$
in $W^{1,2}((0,L_i); \mathbb{M}^{3 \times 3})$ (this follows from trace theorem) and $\whRbb_i^{(j)} (0)=\Rbb_i (0)$. Now take $\Rbb_i^{(j)}= \Pi \whRbb_i^{(j)}$ where $\Pi: \mathbb{M}^{3 \times 3} \to \mathbb{M}^{3 \times 3}$ is a smooth function in the neighborhood of $\SO(3)$ defining projection from the neighborhood of $\SO(3)$ to $\SO(3)$. We define
$$
\yb_i^{(j)} (x_1):= \yb_i (0)+ \int_0^{x_1} \Rbb^{(j)} (s) \eb_1 ds, \quad \db_k^{i,(j)}=\Rbb^{(j)} (x_1) \eb_k,\ \textrm{for} \ k=2,3.
$$
Then $((\yb^{(j)}_1,\db_1^{2,(j)},\db_1^{3,(j)}),\dots (\yb^{(j)}_n,\db_n^{2,(j)},\db_n^{3,(j)})) \in \mathcal{A}$ and $\yb^{(j)}_i$ in $C^2([0,L_i];\ZR^3)$, $\db_i^{2,(j)},\db_i^{3,(j)} \in C^1([0,L_i]; \ZR^3)$ and we also have that $((\yb^{(j)}_{i})',\db_i^{2,(j)},\db_i^{3,(j)})=\Rbb_i^{(j)}$ is converging to $\Rbb_i=
(\yb_{i}',\db_i^{2},\db_i^{3})$ in $W^{1,2} ((0,L_i); \mathbb{M}^{3 \times 3})$. The functions $\betab_i$ are chosen in the following way. We choose
$\alb_i(x_1,\cdot) \in V_i$ (see Remark \ref{napp}) to be the solution of the minimum problem defining $\Qu_i^2 (\Rbb_i^T(x_1) \Rbb_{i}'(x_1))$ (affine function of $\Rbb_i^T(x_1) \Rbb_{i}'(x_1)$). Now take $\betab_i=\Rbb_i \alb_i$  and take $\betab_i^{(j)} \in C^1([0,L_i] \times S_i;\ZR^3)$ defined by convolution (first by first variable and then by last two variables) such that $\betab_i^{(j)} \to \betab_i$ and $\partial_k \betab^{(j)}_{i} \to \partial_k \betab_{i}$ (for k=2,3) in $L^2(\Omega;\ZR^3)$. By an application of the Nemytsky operators theory (see \cite[p.15]{Ambrosseti}) we have that for every $i$
$$
\int_{[0,L_i]} \Qu_i^3 ((\Rbb_i^{(j)})^T (x_2 (\db_{i}^2)'+x_3 (\db_{i}^3)' \ | \ \partial_2 \betab_{i}^{(j)} \ | \ \partial_3 \betab_{i}^{(j)})) dx \to \int_{[0,L_i]} \Qu_i^3 ((x_2 \Rbb_i^T (\db_{i}^2)'+x_3 \Rbb_i^T (\db_{i}^3)' \ | \ \partial_2 \alb_{i} \ | \ \partial_3 \alb_{i}))dx.
$$
Therefore, we can assume (by taking a subsequence) that
$$
\Big| I(\yb,\db^2,\db^3)-
\frac{1}{2}\int_{[0,L_i]} \Qu_i^3 ((\Rbb_i^{(j)})^T (x_2 (\db_{i}^{2,(j)})'+x_3 (\db_{i}^{3,(j)})' \ | \ \partial_2 \betab_{i}^{(j)} \ | \ \partial_3 \betab_{i}^{(j)})) dx \Big| < \frac{1}{j}.
$$
For a given $j$, from the previous part of the proof, we can find $\yb^{(h_j)} \in W^{1,2}(\Omega^{h_j}; \ZR^3)$ such that
$$
\frac{1}{h_j^2} \sum_{i=1}^n \|\yb^{(h_j)}\circ \Qbb_i- D_i(\yb_i^{(j)},\db^{2,(j)}_i,\db^{3,(j)}_i)\|^2_{W^{1,2}((h_j,L_i) \times h_jS_i;\ZR^3)} <\frac{1}{j}$$
and
$$
\Big|\frac{1}{h_j^4} E^{(h_j)}(\yb^{(h_j)})-\frac{1}{2} \int_{(0,L_i)\times S_i} \Qu_3^i((\Rbb_i^{(j)})^T (x_2 (\db_{i}^{2,(j)})'+x_3 (\db_{i}^{3,(j)})'\ | \ \partial_2 \betab_{i}^{(j)} \ | \  \partial_3 \betab_{i}^{(j)}))\Big| < \frac{1}{j}.
$$
By the triangle inequality we have that $\yb^{(h_j)}$ satisfies (\ref{kvg}) and
$$
\Big|\frac{1}{h_j^4} E^{(h_j)}(\yb^{(h_j)})-  I(\yb,\db^2,\db^3) \Big| \to 0.
$$
The case $(\yb,\db^2,\db^3) \notin \mathcal{A}$ is obvious.
\end{prooof}

\section{Minimization}
\setcounter{equation}{0}

Since we can not formulate the problem of junction on a canonic domain in a simple way we have to adapt techniques of $\Gamma$--convergence and use the asymptotics of the infimizing sequence in the form  (\ref{kvg}).

We suppose that the external body force is given by the density $\fb_r^{(h)} \in L^2(\Omega^{h};\ZR^3)$ and that the external surface force is given by the density $\gb^{(h)}_r \in L^2(\partial \Omega^h;\ZR^3)$ (we assume both are dead loads). As is usual in lower-dimensional modeling the scaling of the surface force densities is different at rod ends and the lateral boundary. Therefore we introduce the notation
$$
\gb^{(h)}_{rl}=\gb^{(h)}_r |_{\partial \Omega^h \backslash \cup_{i=1}^n \Qbb_i (\{L_i\} \times h S_i)},  \qquad \gb^{(h)}_{re}=\gb_r^{(h)}|_{\cup_{i=1}^n \Qbb_i (\{L_i\} \times h S_i)}.
$$
We give the result for Neumann boundary condition on the whole domain, i.e for the pure traction problem.
Therefore we suppose that the resultant of all forces is zero, i.e. $\int_{\Omega^h} \fb_r^{(h)} (x) dx+ \int_{\partial \Omega^h} \gb_r^{(h)}(x) dx=0$, and look for the minimum that satisfies $\fint_{\Qbb_1(\{L_1\} \times h S_1)} \yb^{(h)}(x) dx=0$.
\begin{theorem}\label{zadnji}
For every $h$ we define the functional
$$J^{(h)}(\vb)= \int_{\Omega^h} W( \nabla \vb (x)) dx- \int_{\Omega^h} \fb_r^{(h)}(x) \cdot \vb(x) dx-\int_{\partial \Omega^h} \gb^{(h)}_r (x) \cdot \vb(x) dx$$
in the space
$$
V^h=\{ \vb \in W^{1,2} (\Omega^h;\ZR^3) | \fint_{\Qbb_1(\{L_1\} \times h S_1)} \vb(x) dx=0 \}.
$$
Let the scaling of loads is as follows
$$
\fb^{(h)} =\frac{\fb_r^{(h)}}{h^2},\quad  \quad \gb^{(h)}_{l}=\frac{\gb^{(h)}_{rl}}{h^3}, \quad \gb^{(h)}_e=\frac{\gb^{(h)}_{re}}{h^2},
$$
where
$\frac{1}{h} \| \fb^{(h)} \|_{L^2(\Omega^h;\ZR^3)}$, $\frac{1}{h} \| \gb^{(h)}_e \|_{L^2(\cup_{i=1}^n \{L_i\} \times h S_i;\ZR^3)}$ and $\frac{1}{h} \|\gb^{(h)}_{l} \|^2_{L^2(\partial \Omega^h \backslash \cup_{i=1}^n \{L_i\} \times h S_i;\ZR^3)}$ are bounded. Moreover, let us suppose that
\begin{eqnarray}
\label{prvo11} \int_{\Omega^h} \fb_r^{(h)} (x) dx+ \int_{\partial \Omega^h} \gb_r^{(h)}(x) dx&=&0,\\
\label{drugo12} \lim_{h \to 0} \sum_{i=1}^n \frac{1}{h^2} \int_{C_i^{h}} \| \fb^{(h)}(x)-\fb_i ((\Pb^{(h)})^{-1} (\Qbb_i^T x))\|^2 dx&=&0, \\
 \label{trece13} \lim_{h \to 0} \sum_{i=1}^n \frac{1}{h} \int_{\Qbb_i ((h,L_i) \times h \partial S_i)} \| \gb^{(h)}_{l}(x)-\gb_{li}((\Pb^{(h)})^{-1}(\Qbb_i^T x))\|^2 dx&=&0,\\
 \label{cetvrto14} \lim_{h \to 0} \sum_{i=1}^n \frac{1}{h^2} \int_{ \Qbb_i (\{L_i\} \times h S_i)} \| \gb^{(h)}_{e}(x)-\gb_{ei}((\Pb^{(h)})^{-1}(\Qbb_i^T x))\|^2 dx&=&0,
\end{eqnarray}
where  $\fb_i \in L^2((0,L_i) \times S_i; \ZR^3)$, $\gb_{li} \in L^2 ((0,L_i) \times \partial S_i; \ZR^3)$, $\gb_{ei} \in L^2(\{L_i\} \times S_i ; \ZR^3)$, for $i=1, \cdots,n$.

Then we have that
$|\inf_{\vb \in V^h} J^{(h)} (\vb)| \leq Ch^4 $.
Let us take the sequence $\yb^{(h)} \in V^{h}$ that satisfies
\begin{equation} \label{uvjet} J^{(h)} (\yb^{(h)}) \leq \inf_{\vb \in V^{h}} J^{(h)} (\vb)+o(h^4), \end{equation}
($o(h^4)$ means that $\lim_{h \to 0} \frac{o(h^4)}{h^4}=0$).
Let the sequence $(h_j)$ converge to $0$.
Then there exists a subsequence of $(h_j)$ (still denoted by $h_j$) and  $(\yb, \db^2, \db^3) \in \cA$
such that
\begin{equation}\label{kvgista}
\lim_{j \to \infty} \frac{1}{h_j^2} \sum_{i=1}^n \|\yb^{(h_j)}\circ \Qbb_i- D_i(\yb_i,\db^2_i,\db^3_i)\|^2_{W^{1,2}((h_j,L_i) \times h_jS_i;\ZR^3)}=0.
\end{equation}
The limit $(\yb, \db^2, \db^3)$  minimizes the functional
\begin{eqnarray*}  J(\yb, \db^2, \db^3)&=&I(\yb, \db^2, \db^3)-\sum_{i=1}^n \int_{0}^{L_i} \int_{S_i} \fb_i (x) dx_2 dx_3 \cdot \yb_i(x_1) dx_1\\
&& -\sum_{i=1}^n\int_{0}^{L_i} \int_{\partial S_i} \gb_{li}(x) ds \cdot \yb_i(x_1) dx_1 -\sum_{i=1}^n \int_{\{L_i\}\times  S_i} \gb_{ei} (L_i,x_2,x_3)dx_2 dx_3 \cdot \yb_i(L_i)
\end{eqnarray*}
in the space
$V_l=\{ (\yb,\db^2,\db^3)  \in \mathcal{A}: \yb_1(L_1)=0 \}.$
Moreover, the energies converge to the energy of the limit
$$
\lim_{h \to 0} \frac{1}{h^4} J^{(h)}(\yb^{(h)})=J(\yb, \db^2, \db^3).
$$
\end{theorem}
\begin{prooof}
STEP 1 (a priori estimate for the total energy and $\yb^{(h)}$):
Let us estimate $|J^{(h)}(\Id+\ab^{(h)})|$, where $\ib$ is identity mapping and $\ab^{(h)} \in \ZR^3$ is chosen such that $\Id+\ab^{(h)} \in V^h$ (such $\ab^{(h)}$ exists and is unique). Using (\ref{prvo11}) and $W(\Ibb)=0$ we obtain
 \begin{eqnarray*}
|J^{(h)}(\Id+\ab^{(h)})| &=&|\int_{\Omega^h} \fb^{(h)}_r(x) \cdot \Id(x) dx+ \int_{\partial \Omega^h} \gb^{(h)}_r (x) \cdot \Id(x) dx|\\
&\leq&  Ch^3\|\fb^{(h)}\|_{L^2(\Omega^h)}+Ch^{7/2} \|\gb^{(h)}_{l}\|_{L^2(\partial \Omega^h \backslash (\cup_{i=1}^n \Qbb_i(\{L_i\} \times h S_i ))}\\
&&+Ch^3 \|\gb^{(h)}_{e}\|_{L^2(\cup_{i=1}^n \Qbb_i(\{L_i\} \times h S_i ))}\leq Ch^4
 \end{eqnarray*}
Then from (\ref{uvjet}) we conclude that $\frac{1}{h^4} J^{(h)}(\yb^{(h)}) \leq  C $. From this we want to conclude that $\frac{1}{h^4} \int_{\Omega^h} \dist^2 (\nabla \yb^{(h)}, \SO(3))^2 dx< \infty$, so we have to estimate the energy from below
\begin{eqnarray}\nonumber
\frac{1}{h^4} J^{(h)} (\yb^{(h)}) &\geq& C_W \frac{1}{h^4}\int_{\Omega^h} \dist^2( \nabla \yb^{(h)}, \SO(3))\\
\nonumber
&&-\frac{1}{h^2}\|\fb^{(h)}\|_{L^2(\Omega^{h})}\|\yb^{(h)}\|_{L^2(\Omega^{h})}-\frac{1}{h}\| \gb^{(h)}_{l}\|_{L^2(\partial \Omega^h \backslash (\cup_{i=1}^n \{L_i\} \times h S_i )}\| \yb^{(h)} \|_{L^2(\partial \Omega^h \backslash (\cup_{i=1}^n \Qbb_i(\{L_i\} \times h S_i) )}\\
\label{osnovno4}
&& -\frac{1}{h^2}\| \gb^{(h)}_e \|_{L^2(\cup_{i=1}^n \Qbb_i(\{L_i\} \times h S_i) )}\| \yb^{(h)}\|_{L^2(\cup_{i=1}^n \Qbb_i(\{L_i\} \times h S_i ))} \\
\nonumber
&\geq& C_W \frac{1}{h^4}\int_{\Omega^h} \dist^2( \nabla \yb^{(h)}, \SO(3))- C \big(\frac{1}{h} \|\yb^{(h)}\|_{L^2(\Omega^{h})}+\frac{1}{h^{1/2}} \|\yb^{(h)}\|_{L^2(\partial \Omega^h \backslash (\cup_{i=1}^n \Qbb_i(\{L_i\} \times h S_i ))}\\
\nonumber
&&+\frac{1}{h}\| \yb^{(h)}\|_{L^2(\cup_{i=1}^n \Qbb_i(\{L_i\} \times h S_i) )}\big).
\end{eqnarray}
We have
$$ \|\yb^{(h)} \|_{L^2(\Omega^{h})}= \sum_{i=1}^n \|\yb^{(h)} \|_{L^2(C_i^{h})}+\| \yb^{(h)} \|_{L^2(T^{h})} $$
In the same way as in Lemma \ref{lemasranje} we conclude that there exists constant $C$ independent of $h$ (using rescaling $\alpha^{(h)} (x_1,x_2,x_3)=(x_1,h x_2, h x_3)$) such that for every $i$ and every $\yb \in W^{1,2}(\Omega^h;\ZR^3) $ we have
 \begin{eqnarray} \label{josjednom} \left\|\yb-\fint_{\Qbb_i(\{L_i\} \times h S_i)} \!\!\yb \, dx\right\|_{L^{2}(C_i^{h};\ZR^3)}  &\leq& C \|\nabla \yb \|_{L^2(C_i^{h};\ZR^3)}, \\
 \label{josjosjednom} \left\|\yb-\fint_{\Qbb_i(\{h\} \times h S_i)} \!\!\yb \, dx\right\|_{L^{2}(C_i^{h};\ZR^3)}  &\leq& C \|\nabla \yb \|_{L^2(C_i^{h};\ZR^3)}.
  \end{eqnarray}
>From this we conclude that there exists a constant $C$ independent of $h$ such that
\begin{equation} \label{josmalox}
\left\|\fint_{\Qbb_i(\{h\} \times h S_i)} \!\!\yb \, dx-\fint_{\Qbb_i(\{L_i\} \times h S_i)} \!\!\yb \, dx\right\|
\leq \frac{C}{h} \|\nabla \yb \|_{L^2(C_i^{h};\ZR^3)}
\end{equation}
By using scaling $\alpha^{(h)}(x_1,x_2,x_3)=(h x_1, h x_2, h x_3)$ we conclude that there exists $C$ independent of $h$ such that for every $i,l$
\begin{eqnarray}
\label{gotovo}\left\|\yb-\fint_{\Qbb_i(\{h\}\times hS_i)} \!\! \yb\, dx\right\|_{L^{2}(T^h;\ZR^3)} &\leq& hC \| \nabla \yb \|_{L^{2}(T^h;\ZR^3)}, \\
\label{gotovo1} \left\|\fint_{\Qbb_i(\{h\} \times hS_i)} \!\! \yb \, dx-\fint_{\Qbb_l(\{h\} \times hS_l)} \!\! \yb \, dx\right\| &\leq&  \frac{C}{\sqrt{h} } \| \nabla \yb \|_{L^{2}(T^h;\ZR^3)}.
\end{eqnarray}
Using estimate (\ref{josjednom})  for $i=1$ and the fact that $\yb^{(h)} \in V^{h}$  we conclude that
$$
\| \yb^{(h)} \|_{L^2(C_1^{h};\ZR^3)} \leq C \|\nabla \yb^{(h)} \|_{L^2(C_1^{h};\ZR^3)} \leq C (\|\dist(\nabla \yb^{(h)},\SO(3))\|_{L^2(\Omega^{h};\ZR^3)}+h).
$$
Using estimate (\ref{josmalox}) we conclude that
$$\|\fint_{\Qbb_1 (\{h\} \times h S_1)} \yb^{(h)}\|  \leq \frac{C}{h} \|\nabla \yb^{(h)} \|_{L^2(C_1^{h};\ZR^3)}\leq C (\frac{1}{h}\|\dist(\nabla \yb^{(h)},\SO(3))\|_{L^2(C_1^{h};\ZR^3)}+1).
$$
Using estimate (\ref{gotovo}) we conclude
$$ \| \yb^{(h)} \|_{L^{2}(T^h;\ZR^3)} \leq h^{3/2} \|\fint_{Q_1(\{h\} \times hS_1)} \yb^{(h)} \|+Ch\|\nabla \yb^{(h)} \|_{L^{2}(T^h;\ZR^3)}. $$
Since
\begin{equation}
\|\nabla \yb^{(h)} \|_{L^{2}(T^h;\ZR^3)}\leq C( \|\dist(\nabla \yb^{(h)},\SO(3)\|_{L^2(\Omega^h;\ZR^3)}+h^{3/2}), \label{5.31}
\end{equation}
we conclude
\begin{equation} \label{oceli} \| \yb^{(h)} \|_{L^{2}(T^h;\ZR^3)} \leq C(h^{1/2}\|\dist(\nabla \yb^{(h)},\SO(3)\|_{L^2(\Omega^h;\ZR^3)}+h^{3/2}).
\end{equation}
Using estimate (\ref{josjosjednom}) and  (\ref{gotovo1}) for $l=1$ we conclude that for every $i$
$$
\| \yb^{(h)} \|_{L^2(C_i^h;\ZR^3)} \leq  C (\|\dist(\nabla \yb^{(h)},\SO(3))\|_{L^2(\Omega^{h};\ZR^3)}+h).
$$
Thus we have
\begin{equation} \label{osnovno1}
\| \yb^{(h)} \|_{L^2(\Omega^h;\ZR^3)} \leq  C (\|\dist(\nabla \yb^{(h)},\SO(3))\|_{L^2(\Omega^{h};\ZR^3)}+h).
\end{equation}
In the same way one can analyze traces.
First we start from the trace inequality on the cylinder $C_i=(0,1) \times S_i$. For every $\yb \in W^{1,2}(C_i;\ZR^3)$ we have that there exists constant $C$ such that
\begin{eqnarray}
\| \yb- \fint_{\{1\} \times S_i} \yb \|_{L^2(\partial C_i)} \leq C \| \yb-\fint_{\{1\} \times S_i} \yb  \|_{W^{1,2}(C_i;\ZR^3)} &\leq& C\| \nabla \yb \|_{L^{2}(C_i;\ZR^3)},\\
\| \yb- \fint_{\{0\} \times S_i} \yb \|_{L^2(\partial C_i)} \leq C \| \yb-\fint_{\{0\} \times S_i} \yb \|_{W^{1,2}(C_i;\ZR^3)} &\leq& C\| \nabla \yb \|_{L^{2}(C_i;\ZR^3)}.
\end{eqnarray}
By using appropriate scaling and rotation we have that there exists constant $C$ such that for every $C_i^h$  and $\yb \in W^{1,2}(C_i^h;\ZR^3)$ we have
\begin{eqnarray}
\label{prvo}\| \yb-\fint_{\Qbb_i (\{h\}\times h S_i) } \yb \|_{L^2(\Qbb_i(\{L_i\} \times h S_i);\ZR^3)} &\leq& C\|\nabla \yb\|_{L^2(C_i^h;\ZR^3)}, \\
\label{drugo}\| \yb-\fint_{\Qbb_i (\{L_i\}\times h S_i) } \yb  \|_{L^2(\Qbb_i (\{L_i\} \times h S_i);\ZR^3)} &\leq& C\|\nabla \yb\|_{L^2(C_i^h;\ZR^3)}, \\
\label{trece}\| \yb-\fint_{\Qbb_i (\{h\}\times h S_i) } \yb \|_{L^2(\Qbb_i (\{h\} \times h S_i);\ZR^3)} &\leq& C\|\nabla \yb\|_{L^2(C_i^h;\ZR^3)}, \\
\label{cetvrto}\| \yb-\fint_{\Qbb_i (\{L_i\}\times h S_i) } \yb \|_{L^2(\Qbb_i (\{h\} \times h S_i);\ZR^3)} &\leq& C\|\nabla \yb\|_{L^2(C_i^h;\ZR^3)}
\end{eqnarray}
and
\begin{eqnarray}
\label{prvo1}\| \yb- \fint_{\Qbb_i (\{h\}\times h S_i) } \yb  \|_{L^2(\Qbb_i ((h,L_i) \times \partial \, h S_i);\ZR^3)} &\leq& C \frac {1}{h^{1/2}}\|\nabla \yb\|_{L^2(C_i^h;\ZR^3)}, \\
\label{drugo1}\| \yb-\fint_{\Qbb_i (\{L_i\}\times h S_i) } \yb \|_{L^2(\Qbb_i ((h,L_i) \times \partial \, h S_i);\ZR^3)} &\leq& C\frac{1}{h^{1/2}}\|\nabla \yb\|_{L^2(C_i^h;\ZR^3)}, \\
\label{trece1}\| \yb-\fint_{\Qbb_i (\{h\}\times h S_i) } \yb \|_{L^2(\Qbb_i ((h,L_i) \times \partial \, h S_i);\ZR^3)} &\leq& C \frac{1}{h^{1/2}}\|\nabla \yb\|_{L^2(C_i^h;\ZR^3)}, \\
\label{cetvrto1} \| \yb-\fint_{\Qbb_i (\{L_i\}\times h S_i) } \yb \|_{L^2(\Qbb_i ((h,L_i) \times \partial \, h S_i);\ZR^3)} &\leq& C \frac{1}{h^{1/2}}\|\nabla \yb\|_{L^2(C_i^h;\ZR^3)}.
\end{eqnarray}
In the same way we conclude
\begin{eqnarray}
\label{prvo2} \| \yb- \fint_{\Qbb_i (\{h\}\times h S_i) } \yb  \|_{L^2( \partial T^h;\ZR^3)} &\leq& C h^{1/2} \|\nabla \yb\|_{L^2(T^h;\ZR^3)}
\end{eqnarray}
Now by using $\yb^{(h)} \in V^h$ we have from (\ref{drugo}) and (\ref{drugo1}) that
\begin{eqnarray}
\| \yb^{(h)} \|_{L^2(\Qbb_1(\{L_1\} \times h S_1);\ZR^3)} &\leq& C\|\nabla \yb^{(h)}\|_{L^2(C_1^h;\ZR^3)} \leq  C (\|\dist(\nabla \yb^{(h)},\SO(3))\|_{L^2(\Omega^{h};\ZR^3)}+h), \\
\| \yb^{(h)}  \|_{L^2(\Qbb_1 ((h,L_1) \times \partial \, h S_1);\ZR^3)} &\leq& C \frac {1}{h^{1/2}}\|\nabla \yb^{(h)}\|_{L^2(C_1^h;\ZR^3)} \leq
C (\frac{1}{h^{1/2}}\|\dist(\nabla \yb^{(h)},\SO(3))\|_{L^2(\Omega^{h};\ZR^3)}+h^{1/2})
\end{eqnarray}
>From (\ref{prvo}) and (\ref{drugo}) we conclude
\begin{equation}
\label{peto2} \| \fint_{\Qbb_1(\{h\}\times h S_1) } \yb^{(h)} \|  \leq  \frac{C}{h} \|\nabla \yb^{(h)} \|_{L^2(C_1^h;\ZR^3)}.
\end{equation}
>From this and (\ref{gotovo1}) we conclude for every $i$
\begin{equation}
\label{sesto}\| \fint_{\Qbb_i (\{h\}\times h S_i) } \yb^{(h)} \| \leq \frac{C}{h} \|\nabla \yb^{(h)} \|_{L^2(C_1^h;\ZR^3)} +\frac{C}{h^{1/2}} \| \nabla \yb^{(h)} \|_{L^2(T^h;\ZR^3)} \leq \frac{C}{h} \|\nabla \yb^{(h)} \|_{L^2(\Omega^h;\ZR^3)}.
\end{equation}
>From (\ref{prvo1}) and (\ref{osnovno1}) we conclude
\begin{eqnarray}\nonumber
\| \yb^{(h)}  \|_{L^2(\Qbb_i ((h,L_i) \times \partial \, h S_i);\ZR^3)} &\leq& C (\frac {1}{h^{1/2}}\|\nabla \yb^{(h)}\|_{L^2(\Omega^h;\ZR^3)}+h^{1/2}\| \fint_{\Qbb_i(\{h\}\times h S_i) } \yb^{(h)} \| ) \\
\label{deveto}
&\leq& \frac{C}{h^{1/2}} \|\nabla \yb^{(h)}\|_{L^2(\Omega^{h};\ZR^3)} \leq C (\frac{1}{h^{1/2}}\|\dist(\nabla \yb^{(h)},\SO(3))\|_{L^2(\Omega^{h};\ZR^3)}+h^{1/2}).
\end{eqnarray}
>From (\ref{prvo2}) for $i=1$, (\ref{peto2}) and (\ref{5.31}) we conclude
\begin{equation}
\label {osmo} \| \yb^{(h)}  \|_{L^2( \partial T^h;\ZR^3)} \leq C  \|\nabla \yb^{(h)}\|_{L^2(T^h;\ZR^3)} \leq C (\|\dist(\nabla \yb^{(h)},\SO(3))\|_{L^2(\Omega^{h};\ZR^3)}+h^{3/2}).
\end{equation}
>From (\ref{prvo}) and (\ref{sesto}) we conclude for every $i$
\begin{equation}
\label {sedmo} \| \yb^{(h)}\|_{L^2(\Qbb_i(\{L_i\} \times h S_i);\ZR^3)} \leq C\|\nabla \yb^{(h)}\|_{L^2(C_i^h;\ZR^3)} \leq  C (\|\dist(\nabla \yb^{(h)},\SO(3))\|_{L^2(\Omega^{h};\ZR^3)}+h).
\end{equation}
>From (\ref{deveto}), (\ref{osmo}) and (\ref{sedmo}) we conclude
\begin{eqnarray}
\label{osnovno2} \|\yb^{(h)}\|_{L^2(\partial \Omega^h \backslash (\cup_{i=1}^n \Qbb_i(\{L_i\} \times h S_i ))} & \leq & C (\frac{1}{h^{1/2}}\|\dist(\nabla \yb^{(h)},\SO(3))\|_{L^2(\Omega^{h};\ZR^3)}+h^{1/2}),\\
\label{osnovno3}
\| \yb^{(h)}\|_{L^2(\cup_{i=1}^n \Qbb_i(\{L_i\} \times h S_i) )} &\leq&   C (\|\dist(\nabla \yb^{(h)},\SO(3))\|_{L^2(\Omega^{h};\ZR^3)}+h).
\end{eqnarray}
By using  (\ref{osnovno4}), (\ref{osnovno1}), (\ref{osnovno2}) and (\ref{osnovno3}) we conclude that there exists $C_2,  C_3$ such that
\begin{equation}
\label{osnovno6}
C_W \frac{1}{h^4}\int_{\Omega^h} \dist^2( \nabla \yb^{(h)}(x), \SO(3))dx- C_2 \left(\frac{1}{h} \left(\int_{\Omega^h} \dist^2( \nabla \yb^{(h)}(x), \SO(3))dx \right)^{1/2}+1 \right) \leq \frac{1}{h^4} J^{(h)} (\yb^{(h)})  \leq C_3
\end{equation}
Using the fact that for $h \leq 1$
$$
\frac{1}{h} \left(\int_{\Omega^h} \dist^2( \nabla \yb^{(h)}(x), \SO(3))dx \right)^{1/2} \leq    \frac{1}{h^2} \left(\int_{\Omega^h} \dist^2( \nabla \yb^{(h)}(x), \SO(3))dx \right)^{1/2} =: \alpha
$$
we conclude from (\ref{osnovno6}) that
$$
C_W \alpha^2 - C_2 \alpha \leq C_3.
$$
which implies that $\alpha^2$ is bounded, i.e, there exists $C>0$ such that
\begin{equation} \label{dokraja}
\frac{1}{h^4}\int_{\Omega^h} \dist^2( \nabla \yb^{(h)}(x), \SO(3))dx \leq C
\end{equation}
which implies that the left hand side of (\ref{osnovno6}) is bounded as well. This implies $$
|\inf_{\vb \in V^h} J^{(h)} (\vb)| \leq Ch^4.
$$

STEP 2 (the convergence proof for $\yb^{(h)}$ and the scaled total energy):
The estimate (\ref{dokraja}) implies that the assumptions of Theorem~\ref{glavni} (compactness theorem) are fulfilled. Therefore the assumptions of Corollary~\ref{ct3.1l4.1} are fulfilled as well (with $C_{L_1}=0$).
Therefore we conclude that for every sequence $h_j$ there exists a subsequence (still denoted by $h_j$) and  $(\yb,\db^2,\db^3) = ((\yb_1,\db^2_1,\db^3_1), \dots, (\yb_n,\db^2_n,\db^3_n)) \in \mathcal{A}$ and
\begin{equation} \label{zakljucakkonv}
    \lim_{j \to \infty} \frac{1}{h_j^2} \sum_{i=1}^n \|\yb^{(h_j)}\circ \Qbb_i- D_i(\yb_i,\db^2_i,\db^3_i)\|^2_{W^{1,2}((h_j,L_i) \times h_jS_i;\ZR^3)}=0.
\end{equation}
>From this convergence it is obvious that $\yb_1(L_1)= \lim_{j \to \infty} \fint_{\Qbb_1(\{L_1\}\times h_j S_1)}  \yb^{(h_j)} (x) dx=0 $. Thus we have proved that $(\yb,\db^2,\db^3) \in V_l$ and what is left to prove is that it minimizes the functional $J$ in $V_l$. We can use the standard argument from $\Gamma$-convergence, although we have variable domains (and can not apply the $\Gamma$--convergence directly). Let $(\yb_a,\db^2_a,\db^3_a) \in V_l$ and $(\yb_a,\db^2_a,\db^3_a) \neq (\yb,\db^2,\db^3)$. We have to prove that $J(\yb, \db^2, \db^3) \leq J(\yb_a, \db^2_a, \db^3_a)$.
>From the liminf inequality from Proposition~\ref{propo} we conclude
\begin{equation} \label{zaminimum}
I(\yb,\db^2,\db^3) \leq \liminf_{j \to \infty}\frac{1}{h_j^4} E^{(h_j)}( \yb^{(h_j)}).
\end{equation}
By using (\ref{oceli}), (\ref{osmo}) and (\ref{dokraja}) we have
\begin{eqnarray} \label{eq1}
 &&   \frac{1}{h_j^4} \int_{T^{h_j}} \ \left| \fb^{(h_j)}_r(x) \cdot \yb^{(h_j)} (x) \right| dx   \leq  \frac{1}{h_j^2} \| \fb^{(h_j)} \|_{L^2(T^{h_j})} \| \yb ^{(h_j)} \|_{L^2(T^{h_j})} \leq  \frac{1}{h_j^2} \| \fb^{(h_j)} \|_{L^2(\Omega^{h_j})} C(h_j^2+h_j^{3/2}) \to 0 \\
 \label{eq2}
&& \frac{1}{h_j^4} \int_{ \partial T^{h_j} \backslash \cup_{i=1}^n \{ h_j \} \times h_j S_i} \ \left| \gb^{(h_j)}_r(x) \cdot \yb^{(h_j)} (x) \right| dx\\
&& \hspace{10ex}\leq \frac{1}{h_j} \| \gb^{(h_j)}_l \|_{L^2(\partial T^{h_j} \backslash \cup_{i=1}^n \{ h_j \} \times h_j S_i)} \| \yb^{(h_j)} \|_{L^2(\partial T^{h_j} \backslash \cup_{i=1}^n \{ h_j \} \times h_j S_i)}  \leq C(h_j^{(5/2)}+h_j^2) \to 0.
\nonumber
\end{eqnarray}
>From this and  (\ref{drugo12})-(\ref{cetvrto14}), (\ref{osnovno1}), (\ref{osnovno2}), (\ref{osnovno3}), (\ref{zakljucakkonv}) and (\ref{zaminimum}) we conclude that
\begin{equation} \label{konacno11}
J(\yb,\db^2,\db^3) \leq \liminf_{j \to \infty}\frac{1}{h_j^4} J^{(h_j)}( \yb^{(h_j)}).
\end{equation}

Let us by $limsup$-inequality from Proposition~\ref{propo} choose $\yb^{(h_j)}_a$ such that
\begin{equation} \label{konkov} \lim_{j \to \infty} \frac{1}{h_j^2} \sum_{i=1}^n \|\yb_a^{(h_j)}\circ \Qbb_i- D_i(\yb_{i,a} ,\db^2_{i,a} ,\db^3_{i,a} )\|^2_{W^{1,2}((h_j,L_i) \times h_jS_i;\ZR^3)}=0
\end{equation}
and
\begin{equation} \label{svojstvo1} \lim_{j \to \infty} \frac{1}{h_j^4} E^{(h_j)}(\yb_a^{(h_j)})=I(\yb_a,\db_a^2,\db_a^3).
\end{equation}
Let us choose $\ccb^{(h_j)} \in \ZR^3$ such that $\zb_a^{(h_j)}=\yb_a^{(h_j)}+ \ccb^{(h_j)} \in V^{h_j}$. From the convergence (\ref{konkov}) we conclude that $\lim_{j \to \infty} \ccb^{(h_j)} =0$. Thus we have that (\ref{konkov})  is also satisfied for the sequence $\zb_a^{(h_j)}$.
We also see that (\ref{svojstvo1}) is also satisfied for $\zb_a^{(h_j)}$. Therefore, using the lower bound on $W$, it follows that there exists a constant $C$ such that $\sup_{j} \frac{1}{h_j^4} \int_{\Omega^{h_j}} \dist^2 (\nabla \zb_a ^{(h_j)}(x), \SO(3)) dx <C$. In the same way as before we conclude
\begin{equation} \label{konacno12}
J(\yb_a, \db^2_a, \db^3_a)= \lim_{j \to \infty} \frac{1}{h_j^4} J^{(h_j)}(\yb_a^{(h_j)}).
\end{equation}
Finally from (\ref{uvjet}), (\ref{konacno11}) and (\ref{konacno12}) we have
$$ J(\yb, \db^2, \db^3) \leq \liminf_{j \to \infty} \frac{1}{h_j^4} J^{(h_j)} (\yb^{(h_j)}) \leq  \liminf_{j \to \infty} \frac{1}{h_j^4} J^{(h_j)} (\yb_a^{(h_j)})=J(\yb_a,\db^2_a,\db^3_a).$$
That the energies converge can be easily seen by standard argument in $\Gamma$-convergence (we first take sequence $\lb^{(h_j)}$ such that $\frac{1}{h_j^4} J^{(h_j)} (\lb^{(h_j)}) \to J(\yb, \db^2, \db^3)$  and then by using (\ref{uvjet}) conclude that $\lim_{j \to \infty} \frac{1}{h_j^4} J^{(h_j)} (\lb^{(h_j)})=\lim_{j \to \infty} \frac{1}{h_j^4} J^{(h_j)} (\yb^{(h_j)})$. Since this can be done for arbitrary sequence, we have the claim).
\end{prooof}

\begin{remark}\label{eqsile}
Using (\ref{eq1}) and (\ref{eq2}) by straightforward calculation from (\ref{prvo11}) we obtain
$$
\sum_{i=1}^n \left( \int_0^{L_i} \left( \int_{S_i} \fb_i (x) dx_2 dx_3 + \int_{\partial S_i} \gb_{li} (x) ds\right) + \int_{S_i} \gb_{ei} (L_i, x_2, x_3) dx_2 dx_3\right)=0.
$$
This means that in the limit model the total force is zero as well.
\end{remark}

\begin{remark}\label{fix}
Adding a constant to the solution of a pure traction problem gives a solution again, i.e., the set of solutions is closed under translations. Therefore, we had to control behavior of this constant in three-dimensional problem in order to obtain the limit. We did it by requesting that the mean deformation at the end of the first rod (indexed by 1) vanishes. As expected, this constraint results in the limit model in the constraint that the end of the first rod is fixed in the origin ($\yb_1(L_1)=0$). In the limit model we can also consider this constraint as the one which just fixes the translation since again the set of solutions of the pure traction problem is closed under translations.
%
%
%
\end{remark}

\section{Differential formulation of the model}
\setcounter{equation}{0}

In this section we formulate the weak and the differential formulation of the model. It enables us to give the interpretation of the limit model as the model of one--dimensional rods with the transmission conditions at the junction point, see (\ref{p1})--(\ref{Ryoni1}).

Let us define
\begin{eqnarray*}
\ftb_i (x_1) &=&  \int_{S_i} \fb_i (x) dx_2 dx_3 +\int_{\partial S_i} \gb_{li} (x) ds,\\
\Ftb_i  &=&  \int_{S_i} \gb_{ei}(L_i, x_2, x_3) dx_2 dx_3.
\end{eqnarray*}
Then the total energy of the limit model is given by
\begin{eqnarray*}
J(\yb, \db^2, \db^3) = \sum_{i=1}^n \left( \frac{1}{2} \int_0^{L_i} \Qu^i_2 (\Rbb^T_i \Rbb_i') dx_1  - \int_0^{L_i} \ftb_i \cdot \yb_i dx_1 - \Ftb_i \cdot \yb(L_i)\right),
\end{eqnarray*}
for $(\yb,\db^2, \db^3) \in \cA$ and $+\infty$ otherwise.

First, performing partial integration in the force terms in the total energy functional, similarly as in \cite{Muller2}, we remove appearance of $\yb_i$ from the energy functional. In order to do that let us denote
\begin{equation}\label{p}
\ptb_i(x_1) = \int_{x_1}^{L_i} \ftb_i(z) dz + \Ftb_i, \qquad i=1, \ldots, n
\end{equation}
and note that the force equilibrium, according to Remark~\ref{eqsile}, can be expressed by
\begin{equation}\label{c1}
\sum_{i=1}^n \ptb_i(0)=0.
\end{equation}
Then
\begin{eqnarray*}
\sum_{i=1}^n \int_0^{L_i} \ptb_i \cdot \yb_i' dx_1 &=& \sum_{i=1}^n  \ptb_i(L_i) \cdot \yb_i(L_i) - \sum_{i=1}^n \ptb_i(0) \cdot \yb_i(0) - \sum_{i=1}^n \int_0^{L_i} \ptb_i' \cdot \yb_i dx_1\\
&=& \sum_{i=1}^n  \Ftb_i \cdot \yb_i(L_i) - \sum_{i=1}^n \ptb_i(0) \cdot \yb_i(0) + \sum_{i=1}^n \int_0^{L_i} \ftb_i \cdot \yb_i dx_1\\
&=& \sum_{i=1}^n  \Ftb_i \cdot \yb_i(L_i)  + \sum_{i=1}^n \int_0^{L_i} \ftb_i \cdot \yb_i dx_1,
\end{eqnarray*}
since in $\cA$ deformations satisfy $\yb_1(0) = \cdots \yb_n(0)$. Thus the total energy functional can be expressed in terms of $\Rbb_i$, $i=1, \ldots, n$ only, by
$$
\tJ(\Rbb) := J(\yb, \db^2, \db^3) =  \sum_{i=1}^n \left( \frac{1}{2} \int_0^{L_i} \Qu^i_2 (\Rbb^T_i \Rbb_i') dx_1  - \int_0^{L_i} \ptb_i \cdot \Rbb_i \eb_1 dx_1\right),
$$
where we have used the notation $\Rbb = (\Rbb_1, \ldots \Rbb_n)$. Thus we can split the minimization of $J$ in two steps. In the first step we minimize $\tJ$ in the space
\begin{eqnarray*}
 && \cR :=\{ \Rbb = (\Rbb_1, \ldots, \Rbb_n) \in W^{1,2}((0,L_1);\SO(3)) \times \dots \times W^{1,2}((0,L_n);\SO(3)) : \\
 && \hspace{20ex} \Rbb_1(0)\Qbb_1^T=\cdots=\Rbb_n(0)\Qbb_n^T
 \}
\end{eqnarray*}
and in the second step we determine deformations $\yb_i$, $i=1, \ldots, n$ from the equations
\begin{equation}\label{y}
\yb_i' = \Rbb_i \eb_1, \quad \yb_i (0) = \yb_0, \qquad i=1,\ldots, n
\end{equation}
where the constant vector $\yb_0 \in \ZR$ is freely determined from an additional constraint (e.g., $\yb_1(L_1)=0$, see Remark~\ref{fix}).

Now we want to find the weak formulation of the problem
$$
\min_{\Rbb \in \cR} \tJ(\Rbb).
$$
First note that $\Rbb^T_i \Rbb_i'$ are a.e. anti-symmetric matrices. Therefore they possess axial vectors $\ssb_i = \ssb_i(\Rbb_i) \in L^2((0,L_i); \ZR^3)$, i.e.,
$$
\Rbb^T_i \Rbb_i' = \Abb_{\ssb_i}, \qquad i=1, \ldots, n,
$$
where the notation $\Abb_\ssb$ stands for the matrix such that $\Abb_{\ssb} x = \ssb \times x$.
Since $\Qu^i_2$ are quadratic forms of the elements of $\Rbb^T_i \Rbb_i'$ there are positive definite matrices $\Hbb_i$ (positive definiteness of the matrices $\Hbb_i$ follows from the fact that the second derivative of $W$ is greater or equal to $0$ and equal to $0$ exactly on antisymmetric matrices)
such that
$$
\Qu^i_2 (\Rbb^T_i \Rbb_i') = \Hbb_i \ssb_i \cdot \ssb_i.
$$
Thus the total energy functional can be written by
$$
\tJ(\Rbb) =  \sum_{i=1}^n \left( \frac{1}{2} \int_0^{L_i} \Hbb_i \ssb_i(\Rbb_i) \cdot \ssb_i(\Rbb_i) dx_1  - \int_0^{L_i} \ptb_i \cdot \Rbb_i \eb_1 dx_1\right).
$$
In order to obtain the weak and differential formulation of the model we need to find the Gatoux derivative of the functional $\tJ$ over $\cR$. Let $\Rbb \in \cR$, $\eps >0$ and $\vb_i \in C^\infty([0,L_i]; \ZR^3)$, $i=1, \ldots, n$. Let us choose a perturbation $\hRbb=(\hRbb_1, \ldots, \hRbb_n) \in \cR$ of  $\Rbb$ in the following form
$$
\hRbb_i = e^{\eps \Abb_{\vb_i}} \Rbb_i, \qquad i=1, \ldots, n.
$$
In order $\hRbb$ to be in $\cR$ one only needs to fulfill that
$$
\hRbb_i(0) \Qbb_i^T = e^{\eps \Abb_{\vb_i(0)}} \Rbb_i(0) \Qbb_i^T
$$
is independent of $i$. Since $\Rbb_i(0) \Qbb_i^T$ is independent of $i$, $\hRbb_i(0) \Qbb_i^T$  is independent of $i$
if and only if
$$
\vb_1 (0) = \cdots = \vb_n(0).
$$
Thus in the sequel we take
$$
\vb \in \cR_t = \{ \vb=(\vb_1, \cdots, \vb_n) \in C^\infty([0,L_1]; \ZR^3)\times \cdots \times C^\infty([0,L_n]; \ZR^3) : \vb_1 (0) = \cdots =\vb_n(0)\}.
$$
Next we need to compute the axial vectors $\hsb_i$ of $(\hRbb_i)^T (\hRbb_i)'$.
\begin{eqnarray*}
(\hRbb_i)^T (\hRbb_i)' &=& \Rbb_i^T e^{-\eps \Abb_{\vb_i}} (e^{\eps \Abb_{\vb_i'}} \Rbb_i + e^{\eps \Abb_{\vb_i}} \Rbb_i')\\
&=& \Rbb_i^T (\Ibb - \eps \Abb_{\vb_i}+ O(\eps^2)) (\eps \Abb_{\vb_i'} + O(\eps^2)) \Rbb_i + \Rbb_i^T \Rbb_i'\\
&=& \Rbb_i^T \Rbb_i' + \eps \Rbb_i^T  \Abb_{\vb_i'}\Rbb_i + O(\eps^2).
\end{eqnarray*}
Since $\Rbb_i^T  \Abb_{\vb_i'}\Rbb_i x = \Rbb_i^T \vb_i \times x$ (as $\Rbb_i(x_1) \in \SO(3)$ a.e.) we obtain
$$
\hsb_i = \ssb_i + \eps \Rbb_i^T \vb_i' + O(\eps^2), \qquad i=1, \ldots, n.
$$
Now, we plug this perturbation into the functional $\tJ$:
$$
\tJ(\hRbb_i) = \sum_{i=1}^n \left( \frac{1}{2} \int_0^{L_i} \Hbb_i (\ssb_i + \eps \Rbb_i^T \vb_i') \cdot (\ssb_i + \eps \Rbb_i^T \vb_i') dx_1  - \int_0^{L_i} \ptb_i \cdot (\Ibb + \eps \Abb_{\vb_i})\Rbb_i \eb_1 dx_1\right) + O(\eps^2).
$$
Thus the stationary point of the functional $\tJ$ satisfies
$$
\sum_{i=1}^n \left(\int_0^{L_i} \Hbb_i \ssb_i \cdot \Rbb_i^T \vb_i' dx_1  - \int_0^{L_i} \ptb_i \cdot  \Abb_{\vb_i}\Rbb_i \eb_1 dx_1\right) = 0, \qquad \vb \in \cR_t,
$$
i.e.,
$$
\sum_{i=1}^n \left( \int_0^{L_i} \Rbb_i \Hbb_i \ssb_i \cdot \vb_i' dx_1 - \int_0^{L_i} \vb_i \cdot \Rbb_i \eb_1 \times \ptb_i dx_1\right) = 0, \qquad \vb \in \cR_t.
$$
Thus by partial integration on every rod we obtain the equation and the boundary condition
\begin{equation}\label{momenti}
(\Rbb_i \Hbb_i \ssb_i)' + \Rbb_i \eb_1 \times \ptb_i = 0, \qquad \Rbb_i(L_i) \Hbb_i \ssb_i (L_i)  = 0.
\end{equation}
Moreover, since $\vb \in \cR_t$,  we obtain just one condition in the junction point
\begin{equation}\label{c2}
\sum_{i=1}^n \Rbb_i(0) \Hbb_i \ssb_i (0)  = 0.
\end{equation}
Let us now denote
$$
\stb_i = \Rbb_i \ssb_i, \qquad \qtb_i = \Rbb_i \Hbb_i \Rbb_i^T \stb_i.
$$
Then the problem given by (\ref{p}), (\ref{y}), (\ref{momenti}), (\ref{c1}) and (\ref{c2}) can be formulated by
\begin{eqnarray}
\label{p1}
&&\ptb_i' + \ftb_i =0, \quad \ptb_i(L_i) = \Ftb_i, \qquad i=1,\ldots,n,\\
\label{q1}
&&\qtb_i' + \Rbb_i\eb_1 \times \ptb_i = 0, \quad \qtb_i(L_i) =0, \qquad i=1, \ldots, n,\\
\label{z1}
&&\qtb_i = \Rbb_i \Hbb_i \Rbb_i^T \stb_i,\qquad  i=1, \ldots, n,\\
\label{s1}
&&\Rbb_i' =  \Abb_{\stb_i} \Rbb_i, \qquad i=1, \ldots, n,\\
\label{y1}
&&\yb_i' = \Rbb_i \eb_1, \qquad i=1, \ldots, n,\\
\label{pqovi1}
&& \sum_{i=1}^n \ptb_i(0)=0,\qquad \sum_{i=1}^n \qtb_i(0)=0,\\
\label{Ryoni1} &&  \Rbb_1(0) \Qbb_1^T=\cdots = \Rbb_n(0) \Qbb_n^T,\qquad  \yb_1(0) =\cdots = \yb_n(0).
\end{eqnarray}
The first five equations (\ref{p1})--(\ref{y1}) are the equilibrium equations of the nonlinear inextensible rod model, see \cite{Muller1} for the derivation of the model from the threedimensional nonlinear elasticity and \cite{Antman} for the direct foundation of the theory of nonlinear rods. See also \cite{JT} or \cite{Tambaca} for the rod model obtained by linearization of the present one. The model is written as a first order system of ODE's. Introduced unknowns $\ptb_i$ and $\qtb_i$ are the contact force and contact couple corresponding to the $i$-th rod.  The equations (\ref{p1}) and (\ref{q1}) are equilibrium equations together with the boundary conditions. (\ref{z1}) is the constitutional law, (\ref{s1}) and (\ref{y1}) are material restrictions of unshearability and inextensibility. The conditions (\ref{pqovi1})--(\ref{Ryoni1}) are conditions at the junction. The two conditions in (\ref{pqovi1})  are the equilibrium conditions and say that the sum of all contact forces and couples in the junction are 0. The conditions in (\ref{Ryoni1}) are continuity conditions. The first one say that the rotation of the cross-section in the junction is the same looking from all rods. Note here the difference between $\Rbb_i$ and $\Rbb_i \Qbb_i^T$. The matrix $\Rbb_i(0)$ gives actual position of the tangent vector $\Rbb_i(0) \eb_1$ and the cross-section (spanned by $\Rbb_i(0) \eb_2, \Rbb_i(0) \eb_3$).
$\Rbb_i(0)\Qbb_i^T$ is the rotation of the cross-section "in the junction" of the rod for the $i$-th rod (the "difference" between undeformed $\Qbb_i$ and deformed $\Rbb_i(0)$ configuration). The second equation in (\ref{Ryoni1}) say that the deformation in the junction point is the same for all rods.

Thus we conclude that junction (transmission) conditions for junction of rods are given by the equilibrium of contact forces and couples as well as by continuity of the deformations and rotation of the junction.

\begin{remark}
The minimization problem for the total energy $J$ on $\cA$ from Theorem~\ref{zadnji} has at least one solution. Thus $\yb_i \in W^{2,2}((0,L_i); \ZR^3)$, $\Rbb_i \in W^{1,2}((0,L_i); \SO(3))$. From the differential formulation for each rod we can conclude a certain regularity result.
For $\ftb_i \in L^2((0,L_i);\ZR^3)$ one has that $\ptb_i \in W^{1,2}((0,L_i);\ZR^3)$. Therefore $\Rbb_i\eb_1 \times \ptb_i \in W^{1,2}((0,L_i);\ZR^3)$ as well, so $\qtb_i \in W^{2,2}((0,L_i);\ZR^3)$. Using (\ref{z1}) this implies $\stb_i \in W^{1,2}((0,L_i);\ZR^3)$. Now using (\ref{s1}) we obtain that $\Rbb_i \in W^{2,2}((0,L_i);\SO(3))$. Now going back to (\ref{z1}) we obtain that $\stb_i \in W^{2,2}((0,L_i);\ZR^3)$, which again using (\ref{s1}) implies that $\Rbb_i \in W^{3,2}((0,L_i);\SO(3))$ and $\yb_i \in W^{4,2}((0,L_i);\ZR^3)$. This is the most that can be concluded for $L^2$ loads in this fashion.
\end{remark}

\end{document}